\documentclass[11pt]{amsart}
\usepackage{amssymb}
\usepackage{amsmath}
\usepackage{amsfonts}
\usepackage{amscd}
\usepackage{amsthm}
\usepackage{verbatim}
\usepackage[mathscr]{eucal}

\theoremstyle{plain}

 \newtheorem{num}{}[section]
 \newtheorem{thm}[num]{Theorem}
 
 \newtheorem{prop}[num]{Proposition}
 
 \newtheorem{slem}{Lemma}[num]
 \newtheorem{sthm}[slem]{Theorem}
 \newtheorem{sprop}[slem]{Proposition}
 \newtheorem{scor}[slem]{Corollary}

\theoremstyle{definition}
 \newtheorem{ind}[num]{}
 
 \newtheorem{sind}[slem]{}
 \newtheorem{srmk}[slem]{Remark}
 \newtheorem{srmks}[slem]{Remarks}
 \newtheorem{defn}[num]{Definition}

 \numberwithin{equation}{num}

\def\>{\mkern 1mu}
\def\<{\mkern -1mu}
\def\ff#1{\textup{ff}_{\textup{#1}}}
\newcommand{\under}[2]
{\vbox to 0pt{\vskip-#1 ex\hbox{$\scriptstyle #2$}\vss}}
\newcommand{\smcirc}%
  {{\raise.15ex\hbox to.7em{$\hss \scriptstyle\circ\hss$}}} 
 \newcommand{\set}{\!:=}
\newcommand{\iso}%
{{\mkern8mu\longrightarrow \mkern-25.5mu{}^\sim\mkern17mu}}
\newcommand{\osi}%
{{\mkern8mu\longleftarrow \mkern-24.5mu{}^\sim\mkern16mu}}

\def\BA{\, \overline{\!A}\>}

\def\fC{\mathfrak{C}}

\def\C{\mathbb{C}}
\def\cF{\mathcal{F}}

\def\bG{\mathsf{G}}
\def\sG{\mathcal G}
\def\tH{\textup{H}}

\def\sH{\mathcal{H}}
\def\fh{\mathfrak{h}}

\def\cI{\mathcal I}

\def\cL{\mathcal{L}}

\def\fM{\mathfrak{M}}

\def\fm{\mathfrak{m}}

\def\sO{\mathcal{O}}
\def\fp{\mathfrak{p}}
\def\fP{\mathfrak{P}}

\def\fq{\mathfrak{q}}
\def\fQ{\mathfrak{Q}}

\def\OX{{\sO_{\<\<X}}}

\def\dt{\delta}

\def\BX{\>\> \overline{\<\<X\<}\>\>}

\def\B0A{\BA _0}

\def\trdeg{\textup{{tr.deg}}}

\def\spec{\textup{\textsf{Spec}\,}}
\def\specan{\textup{\textsf{Specan}\,}}
\def\proj{\textup{\textsf{Proj}\,}}
\def\pn{\textup{\textsf{Projan}\,}}

\def\ker{\textup{\textsf{Ker}\,}}
\def\tor{\textup{\textsf{Tor}}\>}

\def\triple{f\colon X\to Y}

\def\cond{\textup{($\clubsuit$)}}
\def\condan{\cond_{\textup{an}}}
\def\scond{\textup{($\spadesuit$)}}

\def\ds{\displaystyle}

\def\alb{{\alpha_{\scriptscriptstyle\bullet}}}

\newcommand{\dirlm}[1]%
  {
     {\lim\hskip-1.58em\lower.65ex
       \hbox{$
                {}_{\stackrel{\lower1ex\hbox
                                        {$\scriptstyle -\!\!\!\longrightarrow$}
                                      }{\vbox to0pt{\vss\vskip.6ex
                                            \hbox{$\scriptstyle{}^{#1}$}\vss}}
                   }
            $}
     }
\:}

\begin{document}

\begin{abstract} We investigate conditions for  \emph{simultaneous normalizability}
of a family of reduced schemes, i.e., the normalization of the 
total space normalizes, fiber by
fiber, each member of the family. The main result  
(under more general conditions) is that a flat family  of reduced equidimensional projective $\mathbb C$-varieties $(X_y)_{y\in Y}$ with 
normal parameter space~$Y$---algebraic or analytic---admits a simultaneous normalization if and only if the Hilbert polynomial 
of the integral closure $\>\overline{\<\sO_{\<\<X_{\<y}}\<}\>$ is locally independent
of~$y$. When the~$X_y$ are curves projectivity is not needed, 
and the statement reduces to the well known $\delta$-constant criterion of Teissier.
Proofs are basically algebraic,  analytic results being related via standard techniques (Stein compacta, etc.)~to more abstract  algebraic ones. 
\end{abstract}

\title[Simultaneous Normalization]
{ A Numerical Criterion for\\ Simultaneous Normalization }

\author[H.-J.\:Chiang-Hsieh]{Hung-Jen Chiang-Hsieh}
\address{Dept.\ of  Mathematics, National Chung Cheng University, \\
Chia-Yi 621, Taiwan R.O.C.}
\email{hchiang@math.ccu.edu.tw}
\thanks{H.-J.\:Chiang-Hsieh partially supported by the
National Center for Theoretical Sciences (Taiwan).} 

\author[J.\:Lipman]{Joseph Lipman}
\address{Dept.\ of Mathematics, Purdue University\\
              W. Lafayette IN 47907, USA}
\email {lipman@math.purdue.edu}
\urladdr{www.math.purdue.edu/\~{}lipman/}
\thanks{J.\:Lipman partially supported by the National Security
Agency.
\vspace{.6pt}}

\keywords{simultaneous normalization, Hilbert polynomial}

\subjclass[2000]{14B05,  32S15}
\maketitle

%\centerline{\today}

\section*{Introduction}\label{Intro}
By default, rings are commutative and all schemes and rings are noetherian.

Let $f\colon X\to Y$ be a scheme map which is \emph{reduced} 
(flat, with all nonempty fibers geometrically reduced).
A \emph{simultaneous normalization of  $f\>$}\vspace{.4pt} is a finite map 
$\nu\colon Z\to X$ such that $\bar f\!:=f\<\smcirc\nu$
is \emph{normal} (flat, with all nonempty fibers geometrically
normal), 
and such that for 
each $y\in f(X)$ the induced map of fibers $\nu_y\colon \bar f^{-1}y\to f^{-1}y$ is a
normalization map.% 
\footnote
{A finite map $\mu\colon Z'\to Z$ of noetherian schemes is 
\emph{birational} if it induces an
isomorphism of an open dense subscheme of $Z'$  to an open dense
subscheme of $Z$ (a condition preserved by  flat base change); 
and it is a \emph{normalization map} if in addition, $Z$~
is reduced and $Z'$ is normal.} 
For any base change $Y_1\to Y$ these properties pass over to the projections 
$f_1\colon X_1=X\times_Y Y_1\to Y_1$ and $\nu_1\colon Z_1=Z\times_X X_1\to
X_1$ 
(even if $Y_1$ and $X_1$ are not noetherian). 
If $Y$ is normal then any simultaneous normalization of~$f$ is 
itself a normalization map (Theorem~\ref{Snormal}). When a normalization of $X$ is
a simultaneous normalization of~$f$, we~say that $f$ is \emph{equinormalizable}.
\vspace{.6pt}

\enlargethispage*{15pt}

We will consider families (of fibers) 
given by scheme maps $f\colon X\to Y$ 
subjected to the following fairly mild conditions 
(see \S\ref{Prelims} for information about some of the terms used):

\begin{defn}\label{club} 
A scheme map $f\colon X\to Y$ satisfies $\cond$
if the following all hold:

\begin{itemize} 
\item  $f$ is reduced.\vspace{1pt}

\item For every $y\in f(X)$ the local ring $\sO_{Y\!,\>y}$ is normal and has geometrically normal formal fibers.\vspace{1pt}

\item  $X$ is a formally equidimensional (= locally quasi-unmixed) Nagata scheme. In other words (see~\ref{lqunmixed}), $X$~is a universally catenary Nagata scheme,%
\footnote{\kern.5pt That is, for every affine open subscheme $\spec A\subset X$, $A$ is a universally catenary Nagata ring;
or equivalently, $X$ is covered by affine schemes $\spec A_i$ with each $A_i$ a universally catenary Nagata ring,
see~\ref{univjap} and \ref{unicat} below. The reader who so prefers may simply assume---with minor loss in generality---that $f\colon X\to Y$ is a reduced map of \emph{excellent} schemes, with $X$ equidimensional and $Y$ normal (cf.~\ref{condan}).\looseness=-1
}
 and for every closed point $x\in X$ the local ring $\sO_{\<\<X\<\<,\>x}$ is equidimensional.  \vspace{1pt}

\end{itemize}

\end{defn}

Under \cond\   the flatness condition on $\bar f$ in the definition of simultaneous normalization is superfluous: a normalization  $\nu\colon\BX\to X$ is a simultaneous normalization if and only if for each $y\in f(X)\<$, $\bar f^{-1}y$  is geometrically normal (Corollary~\ref{CorSnormal}). 

\vspace{1.5pt}

\emph{This paper is concerned mainly with  numerical conditions on\/ the fibers of a map~ $f$ satisfying $\cond$ that characterize equinormalizability.} \vspace{1.5pt}

The culminating results, Theorem~\ref{main2} and its analytic avatar Theorem~\ref{main2an}, state roughly that when $f\colon X\to Y$ (satisfying $\cond$)  is projective the sought-after condition is that \emph{the Hilbert polynomial\/ $\sH_{y}(\>\>\overline{\<\<\sO_{\<\<X_y}\<\<}\>\>)$
of the normalization of the structure sheaf on the fiber\/ $X_y\ (y\in Y)$ is locally independent of\/ $y$}.
Since flatness\vspace{-1pt} implies that $\sH_y(\sO_{\<\<X_y})$ is locally 
independent of $y$, one can replace  $\sH_{y}(\>\>\overline{\<\<\sO_{\<\<X_y}\<\<}\>\>)$
here by $\>\sH_{y}(\>\>\overline{\<\<\sO_{\<\<X_y}\<\<}\>\>/\sO_{\<\<X_y})$.\vspace{1pt}

When the fibers of $f$ are one-dimensional, the polynomial $\sH_{y}(\>\>\overline{\<\<\sO_{\<\<X_y}\<\<}\>\>/\sO_{\<\<X_y})$
is  a constant, namely, with $\kappa(y)$ the residue field of $\sO_{Y\!,\>\>y}$, it is the sum 
$$
\sum_{x\in X_y}\delta_y(\sO_{\<\<X_y,\>x})\set\sum_{x\in X_y}\dim_{\kappa(y)}
(\>\>\overline{\<\<\sO_{\<\<X_y,\>x}\<\<}\>\>/\sO_{\<\<X_y,\>x}).
$$ 
In this case, the projectivity assumption on $f$ is not needed, 
and the result reduces to the $\delta$-constant criterion of of Teissier and Raynaud \cite[p.\,73, Th\'eor\`eme 1]{T2}.%
\footnote{That criterion is due to Teissier when $Y$ is the spectrum of a discrete valuation ring---or, in the analytic case, an open disc in $\C$ (cf.~\ref{regular} below). On hearing this, Raynaud quickly generalized it to the case of arbitrary normal $Y\<$. (Some insight into the background of Raynaud's argument might be gleaned from the
introduction to \cite{GS}.) Teissier's result is generalized to deformations of possibly nonreduced curves in \cite[Korollar 3.2.1]{BG}.} That result has largely inspired the
present paper; and though there are gaps in our understanding of the details of the argument presented in \emph{loc.\,cit.,} nevertheless in broad outline our proof runs along similar lines.\vspace{2pt}

In summary,  after recalling some Commutative Algebra in \S1, we do the basics
for equinormalization in \S2, give in \S3 some numerical criteria which, together with\
a supplementary condition (automatically satisfied when the parameter space $Y$
is one-dimensional), characterize equinormalizability, then in \S4 use Quot schemes to eliminate the supplementary condition, thereby obtaining, algebraically, the main results. In \S5 these results are used to prove the corresponding ones for  complex spaces. 
\vspace{2pt}

\enlargethispage*{0 pt}
Actually, the projective case and the one-dimensional case are treated separately. It would be better were we able to prove a theorem for projective maps of formal schemes (resp.~ formal complex spaces), which would yield the results for maps $f\colon X\to Y$ of schemes (resp.~complex spaces) after completion along the closed subspace $N\subset X$ consisting of points  which are not normal on their fibers. (Cf.~\ref{simnorlocal}.) Then the assumption on $f$ could have been weakened to projectivity of the restriction of $f$ to $N$, 
a condition satisfied both when $f$ itself is projective and when the fibers of $f$ are reduced curves (indeed---locally---when $\dim N=0$).\looseness=-2  

A second unresolved issue raised by the main results is why the \emph{global} object~
$\sH_y$ should be involved in the characterization of what is really a \emph{local} phenomenon. To wit:
it follows from~\ref{localgzn}, \ref{unicat}, \ref{localnagata}, \ref{lqunmixed} and~\ref{localgeom} that 
if \cond~holds for~$f$ then for all \mbox{$x\in X\<$,} \cond\ ~holds for the map \mbox{$\<f_{\<x}^{}\colon\<\< \spec\sO_{\<\<X\<\<,\>x}\<\to\<\spec\sO_{Y\!,\>f(\<x\<)}$} induced by~$f$; and by a straightforward globalization of~\ref{univjap} the converse holds iff every integral scheme finite over~$X$ has a nonempty open normal subscheme;
also, with $\nu_{\<\<X}^{}\colon\BX\to X$  a normalization, $f\smcirc\nu_{\<\<X}^{}$ is normal iff with $f_{\<x}^{}$ as above and $\nu_{\<\<X\<\<,\>x}^{}$ the normalization  of $\spec\sO_{\<\<X\<\<,\>x}$, $f_{\<x}^{}\smcirc \nu_{\<\<X\<\<,\>x}^{}$ is normal for all $x\in X\<$.
In fact, the existence of a simultaneous normalization depends, ordinarily, only on the completions of the local rings involved, see~\ref{simnorlocal}. 

The question is, \emph{are there local invariants, somehow related to $\sH_y$, that  characterize equinormalizability at a point?}\vspace{1pt} 

(For one-dimensional fibers,  the $\delta$-invariant provides an affirmative answer.)

\medskip
\centerline{*\quad    *\quad      *\quad      *\quad      *}
\smallskip

Placing our results in a broader context we ask, what sort of information about ``singularity\- type" does normalization convey? In slightly less vague terms,
we'd like to know more about the relation of equinormalization to various versions of \mbox{equisingularity.}
Even for curves equinormalization is weaker, an example being the 
equinormalizable---but not equisingular---family  \mbox{$y^2=tx^2+x^3$} parametrized by $t\in \C$.
(Here $\delta_t\equiv1$.) \vspace{1pt}

The question is,  \emph{which equisingularity conditions  imply equinormalizability?}\vspace{1pt}

 For example, suppose  $f\colon X\to Y$ is a map of $\mathbb Q$-schemes that satisfies $\cond$ and admits a \emph{weak simultaneous resolution,}
that is, there exists a proper birational map $\pi\colon\tilde X\to X $ such that $f\pi$ is smooth and
and for 
each $y\in f(X)$ the induced map of fibers $\pi_y\colon \tilde X_y\to X_y$ is birational. As indicated above, the question of equinormalizability is local, so let us assume, for simplicity, that $X$ and $Y$ are affine. 
We have the Stein factorization through the normalization
$\tilde X\to\BX=\spec\tH^0(\tilde X,\sO_{\<\<\tilde X})\xrightarrow{\nu} X$. That $\nu$  be a simultaneous\- normalization of $f$ means, by~\ref{CorSnormal} below, that ``Stein factorization commutes with passage to the fibers," i.e., that  with $\kappa(y)$ the residue field of~$\sO_{Y\!,\>\>y}$,\looseness=-1
$$
\tH^0(\tilde X, \sO_{\<\<\tilde X})\otimes_R \kappa(y)\cong
\tH^0(\tilde X_y, \sO_{\<\<\tilde X_y})\quad(y\in Y).
$$ 
For this to hold \emph{it suffices that}
$\tH^q(\tilde X, \sO_{\<\<\tilde X})$ \emph{be\/ $R$-flat for all\/ $q\ge0$.} For, these $\tH^q$ can be calculated by the $R$-flat \v Cech complex~$C^\bullet$ associated to a finite affine covering of~$\tilde X\<$; and for any $R$-module $M$ and $p>0$, flatness of $\tH^q$ gives  $\tor_p(\tH^q(C^\bullet), M)=0$, so that the well-known spectral sequences \cite[(6.3.2)]{EGIII}  yield the second of the natural isomorphisms\looseness=-1
$$
\tH^q(\tilde X, \sO_{\<\<\tilde X})\otimes_RM\cong 
\tH^q(C^\bullet)\otimes_R M\cong
\tH^q(C^\bullet\otimes_R M)
\qquad(q\ge0)
$$ 
(where the last term  $\cong \tH^0(\tilde X_y, \sO_{\<\<\tilde X_y}\<)$ when
$q=0$ and $M=\kappa(y)$).

This suggests that \emph{in considering the existence of a simultaneous resolution\/~$\pi\colon\tilde X\to X$ as a condition of equisingularity of the family\/ $f\colon X\to Y$  
one may wish to impose the additional condition that all the higher direct images\/ $R^q\pi_*\sO_{\<\<\tilde X}$ are\/ $Y\<\<$-flat.}\vspace{1pt}

In this connection  
Jonathan Wahl has shown us an example of a one-parameter family of integral two-dimensional quasihomogeneous singularities with
an ``equitopological" simultaneous resolution but no simultaneous normalization. In fact the total space of this family is normal, but one of the members is not normal. Here $R^1\pi_*\sO_{\<\<\tilde X}$ is not flat over the base.  For more in this vein,
see \cite[p.\,341, Thm.\,4.6(i)]{Wa}.

\smallskip

Equinormalizability also holds for \emph{equisaturated} 
families of singularities, in characteristic zero,
see \cite[pp.\;1017--1018]{Z2}. The
arguments given there for hypersurfaces apply somewhat more generally.
Consider a power-series ring
$P\set\mathbb C[[x_1,\dots,x_r,  t_1,\dots, t_s]]$ and a finite
module-free $P$-algebra with $A/(t_1,\dots,t_s)A$ reduced.
Let $\>\>\overline{\<\<A}$ be the integral closure of~$A$.
Let $\nu\colon\spec\>\>\overline{\<\<A}\to\spec A$ and 
$f\colon\spec A\to \spec \mathbb C[[t_1,\dots, t_s]]$ be the
obvious maps.
Suppose 
\emph{the discriminant of the\/ $P$-algebra\/~$A$ $(\<$a principal\/ $P$-ideal\/$)$ is 
generated by a power series independent of the} $t_i$.
Then by a simple extension of \cite[p.\,527, Thm.\,5]{Z1} 
(see \cite[p.\,882, Cor.\,1]{Lp}), $\>\>\overline{\<\<A}$ is of the form 
$A_0[[t_1,\dots, t_s]]$,
so the closed fiber $\spec A_0$ of~$f\nu$ is normal, and by ~\ref{Nishimura}, and \ref{CorSnormal} below, $\nu$ is a simultaneous 
normalization of $f$.
 
 \smallskip
What about \emph{Whitney equisingularity} and equinormalization? For
families of plane curve singularities  Whitney equisingularity is
equivalent to equisaturation, so it does imply equinormalizability.
For certain families of nonplanar curve singularities, the implication still 
holds \cite[p.\,27, Thm.\,III.8]{BGG}.
It would be nice if there were some relation in higher-dimensional  situations; but at present we have no information to offer.

\section{Preliminaries.}\label{Prelims}

In this section we go over briefly some subsequently-used notation, terminology, and basic (long-known) results from commutative algebra.\vspace{1pt}

The total quotient ring of a commutative ring $A$ will be denoted by
$K_{\<\<A}$\kern.5pt, and the integral closure of~ $A$ in $K_{\<\<A}$ by
$\,\overline{\!A}$.
For a prime ideal $\fp \in \spec A$, the residue field of
the local ring $A_{\fp}$ (= field of fractions of $A/\fp$) will be denoted by $\kappa(\fp)$. 
A noetherian local ring $A$  
with maximal ideal~$\fm$ will be denoted by $(A,\fm)$. When we need to refer to the residue field $k_{\<\<A}\!:=\kappa(\fm)$, we may also write
$(A,\fm,k_{\<\<A})$. 

\begin{defn}\label{geornr}
Let $k$ be a field.  A noetherian $k$-algebra $\!A$ is
\emph{geometrically reduced\phantom{$\mkern-3.5mu.$}} (resp.~ \emph{geometrically normal}\kern.5pt)
if the ring $A \otimes_k k'$ is reduced (resp.~normal)%
\footnote{
``$\<\<A$  is normal''  means, by definition,  that  $A_\fp$ is an integrally closed domain
for all $\fp\in\spec A$.}  
for every  field extension  \mbox{$k'\supset k$.} (It~suffices that this be so for every finite, purely inseparable $k'\supset k$, see \cite[(6.7.7), (4.6.1), (6.14.1)]{EGA}. Thus if $k$ is perfect or of characteristic 0, then [$A$~reduced]${}\Rightarrow{}$[$A$ geometrically reduced],
and similarly for ``normal.")
\end{defn}

\begin{sind}\label{localgeom} These geometric properties are local: they 
hold for $A$ if and only if they hold for~$A_\fp$ for all $\fp\in\spec A$.\looseness=-1
\end{sind}

\begin{defn}\label{rnrmap}\textup{
A ring homomorphism $\psi\colon A \to B$ 
 is \emph{reduced} (resp.~\emph{normal}\kern.5pt) 
if $\psi$ is~flat and for every $\fp \in \spec A$ such that $\fp B\ne B$ the
corresponding $\kappa(\fp)$-algebra \mbox{$B\otimes_A\kappa(\fp)$} is
geometrically reduced (resp.~geometrically normal).}\looseness=-1
\end{defn}

\begin{sind}
These properties hold for $\psi$ iff they hold for the induced maps 
$\psi_{\>\fP,\>\fp}\colon A_\fp\to B_\fP$ for all $\fP\in\spec B$ and $\fp=\psi^{-1}\fP$.
\end{sind}

\begin{defn}\label{gzn}\textup{
A noetherian ring $A$ satisfies $\ff{red}$
(resp.~$\ff{nor}$) if for
every maximal \mbox{$A$-ideal} $\fm$ the canonical map
$A_{\fm} \to \widehat{A_{\fm}}$ from the local\vspace{.6pt}
ring $A_{\fm}$ to its completion is
reduced (resp.~normal). [Here ``$\>\ff{}\>$'' is meant to
suggest ``formal fibers.'']}
\end{defn}

\begin{sind}\label{localgzn}\textup{
If $A$ has either of  these $\ff{--}^{}$ properties  then so does any  $A$-algebra essentially of finite type (i.e., any localization of a finitely-generated $A$-algebra), see~\cite[(7.4.5), (7.3.8), (7.3.13)]{EGA}.}
\end{sind}

\begin{defn}\label{nagata}\textup{
A \emph{Nagata ring}%  
\footnote{Also known as \emph{pseudo-geometric ring} or 
\emph{universally Japanese ring}.} 
is a noetherian ring $A$ such that for every $\fp \in \spec A$ and every
finite field extension $L\supset \kappa(\fp)$, the integral
closure of $A/\fp$ in $L$ is module-finite over $A/\fp$.}
\end{defn}

\begin{sind}\label{localnagata}\textup{
 If $A$ is a Nagata ring, then so is any $A$-algebra  essentially of finite type (\cite[p.\,131, (36.1) and p.\,132, (36.5)]{Nag}.}
\end{sind}

\begin{sind}\label{reducednagata}\textup{
 If $A$ is a reduced Nagata ring, then $\BA$ is module-finite
 over $A$.}
\end{sind}

\begin{proof}
 Being noetherian, $A$ has only finitely many minimal
 prime ideals, say $ \fp_1,\fp_2,...,\fp_{\<t}\>$; and being reduced, $A$ is
 naturally identifiable with a subring of $A'\!:=\prod_{i=1}^t \<A/\fp_i$.\vspace{1pt}
 Then any regular element in $A$ is outside $ \fp_i$
 for all $i$ and hence is $A'$-regular.  Thus\vspace{.6pt} we may write
 $A\subset K_{\<\<A} \subset K_{\<\<A'}=\prod_{i=1}^t\kappa(\fp_i)$, and conclude that $ \BA \subset \overline{A'}=\prod_{i=1}^t\<\overline{A/\fp_i}$.
 (Actually, $\BA = \overline{A'}$.)
Since\vspace{-.4pt}  $A$ is a Nagata ring therefore for each $i$, 
 $\overline {A/\fp_i}$ is module-finite over
 $A/\fp_i$, hence over~$A$, and consequently  $\BA \subset
 \prod_{i=1}^t\<\overline{A/\fp_i}$ is module-finite over $A$.
\end{proof}

\begin{sind}\label{univjap} $A$ is a Nagata ring iff $A$ satisfies $\ff{red}$ and for every reduced finitely generated $A$-algebra $A'$ the set of normal points  
is open and dense in $\spec A'\<$, this last condition being implied by $\ff{red}$ when $A$ is semi-local. (It suffices that for all module-finite $A$-algebras $\psi\colon A\to A'$ with $A'$ a domain whose fraction field is finite and purely inseparable over 
that of $A/\ker(\psi)$, $\spec A'$ has a 
nonempty normal open subscheme, see \cite[(6.13.2), (7.6.4), and (7.7.2)]{EGA}.)
\end{sind}

\begin{defn}\label{ucatenary}
A noetherian ring $A$ is \emph{catenary} if for any two
prime $A$-ideals $\fp\supset\fq$ it holds that any two saturated chains of prime ideals between $\fp$ and $\fq$ have the same length (\cite[p.\,31]{Ma}). $A$ is \emph{universally catenary} if
every finite-type $A$-algebra~$B$ is catenary.
\end{defn}

\begin{sind}\label{catenary}
 Let $B\supset A$ be  noetherian integral domains with $B$ a finite-type $A$-algebra. We say  ``the dimension formula holds between $A$ and $B\>$" if for all
 $\fP\in \spec B$ and $\fp=\fP \cap A$,
 $$
 \dim B_{\fP}+\trdeg_{\kappa(\fp)}^{}\kappa(\fP) = \dim A_{\fp}+\trdeg_{K_{\<\<A}}^{} K_{\<\<B}.
 $$
A noetherian ring $A$ is universally catenary iff
the dimension formula holds
between $A/\fp$ and $B$ for every $\fp \in \spec A $ and every integral domain $B\supset A/\fp$ of finite type over $A/\fp$, see \cite[p.\,119, Theorem 15.6]{Ma}. 
(It suffices that this be so for every \emph{minimal} $\fq \in \spec A $, see
\cite[p.\,511, (2.8) and p.\,513, Prop.\,2.13]{R1}.) 
\end{sind}

\begin{sind}\label{unicat}
If the noetherian ring $A$ is universally catenary then so is any $A$-algebra essentially of finite type. Conversely, if $A_\fm$ is universally catenary for every maximal $A$-ideal~$\fm$ then
$A$ is universally catenary. (See \cite[Remarque (5.6.3)]{EGA}.)
\end{sind}

\begin{defn}\label{qunmixed}\textup{
A noetherian  ring $A$ is \emph{equidimensional} if $\dim(A/\fp)=\dim(A)<\infty$ for all minimal prime $A$-ideals~$\fp$; and \emph{biequidimensional} if \emph{any} maximal chain of (distinct) prime ideals in $A$ has length $\dim A<\infty$. $A$ is \emph{formally
 equidimensional}  (or ``locally quasi-unmixed") if the $\fp$-adic completion  of $A_\fp$ is
 equidimensional for every $\fp \in \spec A $. (It suffices that this be so for every
 \emph{maximal} $\fp \in \spec A $, see \cite[p.\,251, Thm.\,31.6]{Ma}).}
\end{defn}

\begin{sind}\label{lqunmixed}\textup{
 A noetherian ring $A$ is formally equidimensional
 if and only if $A$ is universally catenary
  and $A_{\fm}$ is equidimensional for each maximal ideal $\fm$.
 (This is the equivalence of (1) and (2) in \cite[p.\,104, Corollary 2.7]{R2}.)}
So if $A$ is formally equidimensional then for every $\fp\in\spec A$ the local ring $A_\fp$, being  formally equidimensional, is equidimensional. It also follows that if $A$ has only one minimal prime ideal, or if $A$  is biequidimensional, then  $A$ is formally equidimensional if and only if $A$ is universally catenary.
\end{sind}

\begin{sind}\label{unmixed}
If $A$ is reduced and equidimensional, and $B$ is a noetherian $A$-algebra integral over~$A$ and torsion-free as an $A$-module, then $B$ is equidimensional. Indeed, torsion-freeness implies that the inverse image $p$ in $A$ of a minimal prime $B$-ideal $P\/$ is a minimal prime $A$-ideal, and $B/P$ is integral over $A/\<p\>$, so
$\dim B/P=\dim A/p=\dim A$.\looseness=-1
\end{sind} 

\begin{sind}\label{finitefed}

Suppose $A$ reduced. If $A$ is formally equidimensional (resp.~biequidimensional and formally equidimensional), then so is any module-finite torsion-free $A$-algebra $B$.
Indeed, $B$ is universally catenary by \ref{unicat}, so it suffices, because of~ 
\ref{lqunmixed}, to show that $B_\fP$ is equidimensional (resp.~equi\-dimensional of 
dimension $\dim A$) for any maximal $B$-ideal $\fP$. Let $\fp$ be the inverse image of
 $\fP$ in~$A$, a maximal $A$-ideal, and let $P\subset\fP$ and $p\subset \fp$ be minimal primes as in~\ref{unmixed}. Since $A$ is universally catenary, the dimension formula holds between $A/p$ and its finite extension ring $B/P$, so that
  $$
  \dim B_{\fP}/PB_{\fP}=\dim(B/P)_{\fP} +
  \trdeg_{\kappa(\fp)}\kappa(\fP)=
  \dim(A/p)_{\fp} + \trdeg_{\kappa(p)}\kappa(P)= \dim A_{\fp}\>,
  $$
whence the  conclusion.
 \end{sind}

  \begin{ind}
The following results (see
  \cite[p.\,157,  (2.4)]{Ni}), deducing properties of general fibers of certain flat maps from those of closed fibers, will be very useful. 
  
\end{ind}
\begin{sthm} \label{Nishimura}
 \textup{(Nishimura, Andr\'e)} If \,$\phi \colon(R,\fm,k_R)\to(A,\fM,k_{\<\<A})$ is a local homomorphism such that\vspace{1pt}
 
\textup{(i)} $R$ satisfies $\ff{red}$ $($resp.~$\ff{nor}),$ 

\textup{(ii)} the $k_R$-algebra $A/\fm A$ is geometrically reduced $($resp.~geometrically normal\/$\kern.5pt),$ and

\textup{(iii)}  $\phi$ is flat, \\[1pt]
 then the map $\phi$ is reduced $($resp.~normal\/$\kern.5pt)$.
\end{sthm}

\section{Equinormalization.}
 The main results of this section, Theorem~\ref{Snormal} and its corollaries, 
characterize equinormalizability of 
a map~$f\colon X\to Y$ satisfying \cond\ (see Introduction) 
in terms of a normalization of $X\<$, and show that equinormalizability is a strictly local property.\looseness=-1
 
\begin{ind} We begin with a few facts relating condition \cond\ and base change. Terminology remains as above. A generic point of irreducible component of a scheme~$W$ is simply called a \emph{generic point of}~$W$.
\end{ind}
 
\begin{sprop} \label{factA}
 Let\/ $\triple$ satisfy \cond. Let\/ $\mu\colon Z\to X$ be a finite birational~map%
  \footnote
 {See Introduction. Birationality also means $\mu$ induces a bijection from generic points~$\zeta$ of $Z$ to generic points of~$X$ such that for each~$\zeta$, $\sO_{\<\<X\<,\>\mu(\zeta)} \to\sO_{Z ,\>\zeta}$ is an isomorphism \cite[p.\,312, (6.6.4)(ii)]{EGA1}.\looseness=-1}
with $Z$ reduced.  
 Let\/ $Y_1$ be a reduced, essentially-finite-type\/ $Y\<\<$-scheme,%
  \footnote{\kern.5pt In other words, $Y$ has a covering by affine open subsets 
  $\spec A_i$ whose inverse image in $Y_1$ is covered by affine open subsets 
  $\spec A'_{ij}$ with each 
  $A'_{ij}$ a localization of a finitely-generated $A_i$-algebra. 
 }
 and for any $Y\<\<$-scheme\/~$W\<$ set\/ 
 $W_1\!:=W\times_Y Y_1$. Then$\mkern1.5mu:$ \vspace{1pt}

{\rm(i)} $X_1$ is a reduced universally catenary Nagata scheme, for which any normalization map\/~$\nu^{}_1$ factors uniquely as\vspace{-4.5pt}
$$
\overline{\<X_1\<\<}\>\> \xrightarrow{\alpha^{}_{Z,Y_1}}Z_1 
  \xrightarrow{\mu_{\<1}\>\set\mu\times 1\,} X_1.
$$

{\rm(ii)} The inverse image  of the nonnormal locus of\/ $X$ is a nowhere dense closed subset of\/ $X_1$.
 \vspace{2pt}
 
{\rm(iii)} If\/ $Y_1$ is normal then\/ the projection 
 $f^{}_{\<1}\colon X_1\to Y_1$ satisfies\/ \cond.\vspace{1pt}

{\rm(iv)} For all\/ $x\in X$ and\/ $y\!:=f(x),$ if\/
 $\hat f_{\<x}^{}\colon\spec
\widehat{\sO_{\<\<X\<\<,\>x}}\to\spec \widehat{\sO_{Y\!,\>\>y}}$ is the map\vspace{.6pt} 
induced by\/ $f$ $($where $\ \widehat{}\;$ denotes completion$)$
then\/ $\hat f_{\<x}^{}$~satisfies\/ \cond.
More generally, \vspace{.4pt}let\/ $\phi\colon R\to A$ be a ring homomorphism such
that the corresponding scheme map\/ 
$f\colon\spec A\to~\spec R$ satisfies \cond.
Let $I$ be an\/ $R$-ideal such that\/ $R/IR$ is a finite-dimensional
Nagata ring, and 
let\/ $J\supset IA$ be an $A$-ideal. Let\/ $\hat R$
$($resp.~$\hat A)$ be the\/ $I$-adic completion of\/ $R$\vspace{.6pt}
$($resp.~$J$-adic completion of~$A)$. Then the map\/ 
$\hat f\colon\spec \hat A\to\spec\hat R$ induced by\/ $f$ satisfies
\cond.
\end{sprop}

\begin{proof} That $X_1$ is reduced is given by \cite[p.\,184, Cor.\,(ii)]{Ma} or by \cite[(6.15.10)]{EGA};%
\footnote{Replace the map $f$ (resp.~$g$) in \emph{loc.\,cit.}~by our $Y_1\to Y$ 
(resp.~our $f\colon X\to Y$); and  note that in the first paragraph of the proof there, $g$ need only be \emph{reduced}.
} 
that $X_1$ is universally catenary, by~\ref {unicat}; and  that $X_1$ is a Nagata scheme, by~\ref{localnagata}.\vspace{1pt}

In particular, $X$ is a reduced Nagata scheme, so its nonnormal locus is closed 
(\ref{univjap}). Hence, to prove (ii) we need only show that \emph{if\/ $x\in X$ is the image of a generic point\/ ~$x^{}_1$ of\/~$X_1$ 
then\/ $X$ is normal at\/ $x$.} 

If $y\set f(x)$, then the genericity of $x^{}_1$ in $X_1$---hence in its fiber over $Y_1$---implies that $x$~is a generic point of
the fiber $f^{-1}y$. Since $f^{-1}y$ is geometrically reduced, one sees that 
the closed fiber of the flat map
$f_{\<x}\colon\spec\sO_{\<\<X\<,\>x}\to\spec \sO_{Y\!,\>\>y}$ induced by~$f$ 
(namely, $\spec$ of the residue field $\kappa(x)$) is geometrically normal.  
So by~\ref{Nishimura} the map $f_{\<x}$ is normal;  and by  \cite[p.\,184, Cor.\,(ii)]{Ma}  $\sO_{\<\<X\<,\>x}$~ is normal, proving (ii).

For the rest of (i), let $V_1=\spec A\subset X_1$ be an affine open subscheme. $A$ is a reduced noetherian ring, and 
$\nu_1^{-1}V_1=\spec\BA$. Since $\mu_1$ is finite, 
$\mu_1^{-1}V_1=\spec B$ with $B$ a module-finite  $A$-algebra. Again,  the normal points of $X$  form an open dense subscheme~$U\<$; and $\mu$ induces an isomorphism from $\mu^{-1}U$ to $U\<$, hence from $\mu_1^{-1}U_1$ to~$U_1$. Furthermore, as we have just seen, $U_1$ contains every generic point of $X_1$. Thus $\mu_1$ is an isomorphism over a neighborhood of any generic point of~$X_1$,  so for every minimal prime $A$-ideal $\fp$ there is a unique minimal prime $B$-ideal~$\fP$
with inverse image $\fp\subset A$, and a natural isomorphism\vspace{.5pt} $A_\fp\cong B_\fP$. In particular,\vspace{-1pt} the natural map $A\to B$ is injective.
Since $A$ is reduced\vspace{-1pt}
 we have
 $K_A=\prod_\fp A_\fp\cong\prod_\fP B_\fP$.
 Let $B\xrightarrow {\alpha\>}K_A$ be the natural map, so that $A\subset \alpha(B)\subset K_A$. \vspace{.5pt}Since $B$ is finite over~$A$, therefore $\alpha(B)\subset\BA$, so we have a  factorization of the normalization map $A\hookrightarrow\BA$. 
 
 \pagebreak[3] 
 
This factorization is unique. Indeed,\vspace{1pt} since  $A_\fp\to B\otimes_A A_\fp$ is an isomorphism for each minimal prime $A$-ideal $\fp$ therefore  the annihilator
 of the $A$-module $B/A$ is contained in no such~$\fp$, and so it contains a regular element 
$h\in A$ such that $h\alpha(B)\subset A$; and if 
$A\to B\xrightarrow{\under{.95}{\alpha'}}\BA$ is\vspace{.6pt} 
another factorization of $A\hookrightarrow\BA$, 
and $b\in B$, then with this $h$ we have 
$$
h\alpha(b)=\alpha(hb)=\alpha'(hb)=h\alpha'(b)
$$
whence, $h$ being a unit in $K_A$,  
$\alpha(b)=\alpha'(b).$

\enlargethispage*{1pt}
Finally, a simple pasting argument  gives the existence and uniqueness of the asserted global factorization of~$\nu^{}_1$. 
\vspace{1pt}

 (iii) From~\ref{localgzn} it follows that all the local rings of points on $Y_1$ satisfy $\ff{nor}$. It is clear that $f^{}_1$ is flat and has geometrically reduced fibers (any such fiber being obtained from a fiber of $f$
 by base change to a finitely-generated field extension). As in (i),
 $X_1$ is a universally catenary Nagata scheme. For the equidimensionality condition (local on~$X_1$), we may assume that  $Y=\spec R\>$ with $R$ a normal domain, and that $Y_1=\spec R_1$ where $R_1$ is a homomorphic image of a localization at a prime ideal of a polynomial ring $R[T_1,T_2,\dots,T_n]$. As \cite[p.\,250, 31.5]{Ma} takes care of the ``homomorphic image" part, we reduce inductively to the case where $Y_1=\spec R_1$ with $R_1$ the localization
at a prime ideal of the one-variable polynomial ring $R[T\>]$.
 
So, in view of~\ref{lqunmixed}, we need only show that if $A$ is a noetherian ring all of whose localizations are equidimensional then so is $B\!:=A[T\>]$---in other words, if $\fP$ is a prime $B$-ideal then $B_\fP$ is equidimensional. We can replace $A$ by its localization at $\fp\!:=\fP\cap A$, so we may assume that $A$ is local, with
maximal ideal $\fp$.  Then $B/\fp B=(A/\fp)[T\>]$, and so either (a): $\fP=\fp B$ or 
(b): $\dim B_\fP/\fp B_\fP=1$. \vspace{1pt}

Let $\fQ\subset\fP$ be a minimal prime $A[T\>]$-ideal, necessarily equal to $\fq A[T\>]$ with
$\fq\!:=\fQ\cap A$  a minimal prime $A$-ideal. Then by \cite[p.\,116, 15.1]{Ma}, 
$\dim B_\fP/\fQ= \dim A/\fq$ in case (a), and  $\dim B_\fP/\fQ= \dim A/\fq+1$
 in case (b). In either case, $B_\fP$ is equidimensional.\vspace{1pt}

(iv) Any maximal $\hat A$-ideal  is of the form $\fm\hat A$ for some maximal 
$A$-ideal $\fm$, and the $\hat \fm$-adic completion of $\hat A$ is the
$\fm$-adic completion of $A$, so is equidimensional. 
Thus  $\hat A$ is formally
equidimensional (see~\ref{qunmixed}), hence universally catenary (see
~\ref{lqunmixed}). Moreover, by \cite[p.\,107, 2.3]{Mr}, 
$\hat A$ is a Nagata ring.

By \cite[(7.4.6)]{EGA} the fibers of the canonical map 
$\spec\hat R\to\spec R$ are geometrically normal, whence by \cite[p.\,184,
Cor.\,(ii)]{Ma}, $\hat R$ is normal. If $\fq\in\spec\hat R$ is the $\hat
f$-image of $\fp\in\spec\hat A$ then $\hat A_\fp$ is a faithfully flat
$\hat R_\fq$-algebra (\cite[p.\,103]{Bo}), whence the polynomial ring
$\hat A_\fp[T\>]$ is a faithfully flat $\hat R_\fq[T\>]$-algebra. As in
the proof of (iii), since $\hat A$ is formally equidimensional, $\hat
A_\fp[T\>]$ is locally equidimensional and universally catenary (see
~\ref{lqunmixed}).  By \cite[p.\,252, 31.7]{Ma} and \cite[p.\,250,
31.5]{Ma}, it follows that $\hat R_\fq$ is universally
catenary. Since $\hat R\>/I\hat R\>\cong R/IR$ is finite dimensional and
$\hat R$ is a catenary normal ring all of whose maximal ideals contain
$I\hat R$,
therefore $\hat R$ is finite dimensional; and moreover the formal fibers of
the Nagata ring $\hat R\>/I\hat R$ (i.e., of $R\>/IR$) are geometrically normal. Now 
\cite[Satz~2]{BR} gives that the formal fibers of $\hat R$ are
geometrically normal.

It remains to show that $\hat f$ is reduced. By \cite[p.\,103]{Bo}, $\hat A$ is a flat $\hat R$-algebra. Whether every fiber $\hat f^{-1}\fq\ (\fq\in\spec\hat R)$ is geometrically reduced\vspace{-.6pt}  is a local question (see~\ref{localgeom}),  so  we need only show for each $\fp$ as above
that $\hat A_{\fp}/\fq\hat A_{\fp}$\vspace{-.6pt} is geometrically reduced. There is a maximal $\hat A$-ideal $\fp'$ containing $\fp$, whose $\hat f$-image is a prime $\hat R$-ideal
$\fq'\supset I\hat R$; and it will suffice then to show that all the fibers of  the map\vspace{-1pt} 
$\spec \hat A_{\fp'}\to\spec\hat R_{\fq'}$ induced by~$\hat f$ are geometrically reduced. In view of ~\ref{Nishimura}, it will even suffice to show that $\hat f^{-1}\{\fq'\}$ is
geometrically reduced.

%\enlargethispage*{0pt}
For this, let  $\fq''$ be the inverse image in  $R$ of $\fq'$. Since $\hat R\>/I\hat R\cong R\>/IR$ and 
$\fq'\supset I\hat R$ therefore $\fq'=\fq''\hat R  +I\hat R$; and as $I\hat R$ is the Jacobson radical of $\hat R$, this means that $\fq'=\fq''\hat R$, whence 
$\hat A/\fq'\<\hat A$~is the $J$-adic completion of $A/\fq''\<A$.

By~(i), $A/\fq''\<A$ is a reduced \mbox{Nagata} ring,
so  $A/\fq''\<A$ satisfies~$\ff{red}$ (see~\ref{univjap}). Consequently, 
the fibers of the~natural map $\spec( \hat A/\fq'\<\hat A)\to\spec(A/\fq''\<A)$ are
geometrically reduced \cite[(7.4.6)]{EGA}, whence $\hat A/\fq'\<\hat A$~is
reduced\vspace{.6pt} \cite[p.\,184, Cor.\,(ii)]{Ma}. Thus $\hat f^{-1}\{\fq'\}$---the $\spec$ of a localization of ~$\hat A/\fq'\<\hat A$---is indeed reduced.

Similar considerations apply if $\hat A/\fq'\<\hat A$ and $A/\fq''\<A$ are
replaced by their respective tensor products with a finite  integral
$\hat R/\fq'(\,\cong R/\fq'')$-algebra whose fraction field is purely inseparable over that of
$R/\fq''$.
It follows that $\hat f^{-1}\{\fq'\}$ is \emph{geometrically} reduced,
as desired.
\end{proof}

\begin{ind}\label{dominant} Here are some conditions under which the map $\alpha\set\alpha^{}_{Z,Y_1}$\vspace{.5pt} in~\ref{factA}  is  \emph{schematically dominant,}  i.e., the associated map\vspace{.6 pt} $\sO^{}_{\<\<Z_1}\to\alpha_{*}\sO_{\overline{\<\<X_1\<\!}}\>\>$ is injective
(\cite[p.\,284 (5.4.2)]{EGA1}). This condition is easily seen to be equivalent to
``$Z_1$  reduced  and $\mu_1$  birational."
\end{ind}

%:Prop 2.2

For any  $y\in Y$ and any scheme map $g\colon W\to Y\<$, $\kappa(y)$
denotes the residue field of the local ring\/ $\sO_{Y,\>\>y}$ and $W_y$
denotes the fiber $g^{-1}y\set W\times_Y\spec \kappa(y)$. Recall that
an \emph{associated point} of a scheme~$W$ is one in whose local ring every
nonunit is a zerodivisor. 

\begin{sprop}\label{birational}  With notation and assumptions as in\/ \textup{\ref{factA},} the following conditions are equivalent.\vspace{1.5pt}

{\rm(i)}  For all\/ $Y_1,$  $Z_1$ is reduced and\/ $\mu_1$ is birational.\vspace{1pt}

{\rm(ii)} For all\/ $y\in Y\<,$ $\mu\times 1$ maps any associated point of\/ $Z_y$ to a generic point of\/~$X_y.$

{\rm(ii)$'$} Condition\/ \textup{(ii)} holds for all\/ $y=f(x)$ with\/ $x$ a closed point of\/~$X\<$.\vspace{1pt}

{\rm(iii)}  For all\/ $z\in Z,$ with\/ $x\set \mu(z)$ and\/ $y\set f(x)$  the local ring\/ $\sO_{\<Z_y,\>z}$ is reduced and equidimensional of dimension
 $\dim\sO_{\<\<X_y,\>\>x\>}$.

{\rm(iii)$'$} For all\/ $y=f(x)$ with\/ $x$  closed in\/~$X\<,$
every associated point of\/ $Z_y$ is generic.\vspace{1pt}
 
{\rm(iv)} The map\/ $f\<\mu\colon Z\to Y$ satisfies\/ $\cond.$\vspace{1pt}

{\rm(iv)$'$} The map\/ $f\<\mu\colon Z\to Y$ is reduced.\vspace{1.5pt}

\noindent When these conditions hold, 
the cokernel of the natural map\/  $\sO_{\<\<X}\to \mu_*\sO_Z$ is $Y\<\<$-\kern.5pt flat.

\end{sprop}

\begin{proof} 
It is clear that (i) for $Y_1=\spec \kappa(y)$ implies (ii), and trivial that
(ii)${}\Rightarrow{}$(ii)$'$, (iii)${}\Rightarrow{}$(iii)$'$ 
and~(iv)${}\Rightarrow{}$(iv)$'\Rightarrow{}$(iii)$'$. 
To establish all the equivalences it suffices then to show that 
(ii)${}\Rightarrow{}$(iii), (iii)$'\Rightarrow{}$(ii)$'$,
(ii)$'\Rightarrow{}$(i),  and (ii)$'\Rightarrow{}$(iv). 
Along the way, in Lemma~\ref{flat}, we will show 
that the last assertion in the Proposition follows from (ii)$'$.
\vspace{1.5pt}

(ii)${}\Rightarrow{}$(iii).
Let $R\set \sO_{Y,\>y\>}$, with maximal ideal $\fm$, and let $A\set \sO_{\<\<X\<,\>x}$.
Since $\mu$ is finite and birational, and $Z$ is reduced, therefore 
$Z\times_X\spec A=\spec B\>$ where $A\subset B\subset \,\overline{\!A}$.
 Set $A_1\set A/\fm A$, $B_1\set B/\fm B$. 
By~\cond, $A_1$ is reduced, i.e., $\fm A$ is an intersection of prime
ideals.
Since $B$ is a finite extension of $A$, it follows that $\fm A=\fm
B\cap A$, i.e., the natural map $A_1\to B_1$ is injective, making
$B_1$ a finite $A_1$-module. Hence  $\dim B_1=\dim A_1$. 

If (ii) holds then $B_1$ is a torsion-free
$A_1$-module, because  any regular element  $\bar h\in A_1$ lies in no minimal
prime of $A_1$, hence in no associated prime of $(0)$ in $B_1$, 
so that $\bar h$ is $B_1$-regular.
That $B_1$ is reduced follows easily from
the existence, given by Lemma~\ref{hexist} below, of 
an $A_1$-regular---hence~$B_1$-regular---element $h\in A$ such that $hB\subset A$ (whence
$hB_1\subset A_1$).  
By \cond\ and~\ref{lqunmixed}, $A$~is universally catenary, equidimensional, and flat over ~$R$; 
so~by \cite[p.\,250, 31.5]{Ma}, $A_1$ is equidimensional,
and it is also universally catenary (since~$A$~is), 
hence biequidimensional. Then \ref{finitefed}
shows that $B_1$ is biequidimensional, and since $\dim B_1=\dim A_1$, 
(iii)  follows. 

\pagebreak[3]

\begin{small}
[Conversely, assuming (iii) we have in the above
 situation that $B_1$ is reduced and equidimensional, and
if $\>\fP$ is a minimal prime of $B_1$ and $\fp\!:=\fP\cap A_1$ then, 
since $B_1$ is an integral extension of $A_1$ and
 $B_1/\fP$ is an integral extension of $A_1/\fp$, 
 $$
\dim B_1=\dim B_1/\fP=\dim A_1/\fp \le \dim A_1=\dim B_1,
$$
whence $\fp$ is a minimal prime of $A_1$; and  (ii) follows.]
\vspace{1pt}
\end{small}

\begin{slem}\label{hexist}
  Under the preceding circumstances, there exists\/ $h \in A$ such that\/$:$ \\
  \ \textup{(a)} $0:_Ah=0$.\\
  \ \textup{(b)} $\fm A:_A h=\fm A$. \\
  \ \textup{(c)} $h \BA \subset A$.\\
  \ \textup{(d)} $h\>\overline{\<A_1\<\<}\>\>\subset A_1$.
\end{slem}

\begin{proof}
Let $\fC=A:_A\BA$ be the conductor of $A$, let $\fC_1$ be the
conductor of~$A_1$, and with $\pi\colon A\to A_1$ the natural map 
let $\fC'\set\pi^{-1}\fC_1$. 
Since both $A$ and $A_1$ are reduced (see~\ref{factA}\kern.5pt (i)), 
it suffices that there exist $h \in \fC\cap\fC'$ such that 
$h$  belongs to no minimal prime divisor $\fP$ of~(0) or of $\fm A$, or equivalently,
that $\fC\cap\fC'$ is not contained in any such~$\fP$, i.e., neither $\fC$ nor $\fC'$ is
contained in any such $\fP$.
  
  If $\>\fP$ is a minimal prime divisor of~$(0)$ 
then the field $A_\fP$ is normal; 
and if $\>\fP$ is a minimal prime divisor of $\fm A$ then $A_\fP$ is
normal, by~\ref{factA}\kern.5pt(ii) for $Y_1=\spec(R/\fm)$. 
In either case, since $\BA$  is module-finite over the Nagata ring $A$ 
(see~\ref{reducednagata}) therefore $\fC_\fP$ is the conductor of
$A_\fP$, so $\fC\not\subset\fP$. \looseness=-1 
  
  As $\overline{A_1\<\<}\>\>$ is module-finite over the 
reduced Nagata ring $A_1$, $\fC'$ cannot be contained in 
any minimal prime divisor of~$\fm A$. Finally, if $\fP$ is a 
minimal prime divisor of~$(0)$ then $\fC'\not\subset \fP$ 
because the natural map $\phi\colon R\to A$ is flat, so
$\phi^{-1}\fP=(0)$ and  
$$
  \fC'\<\subset \fP\implies \fm A\subset \fP\implies \fm=\phi^{-1}(\fm A)
  \subset\phi^{-1}\fP=(0)\implies \<\<A_1=\<A\implies\fC'\<=\fC\not\subset \fP.
$$
\end{proof}

(iii)$'\Rightarrow{}$(ii)$'$.   
It suffices to show that if 
$z$ is a generic point of $Z_y$ then $x'\set \mu(z)$ is a generic point of~$X_y$. Let $(R,\fm)\set\sO_{Y\!,\>\>y}\>$, with maximal ideal~$\fm$, $A'\set\sO_{\<\<X\<,\>x'}$, 
$B'\set\sO_{Z,\>z}$. It follows from~\ref{finitefed} that $\dim B'=\dim A'$.
Then \cite[p.\,116, 15.1]{Ma} gives
$$
0=\dim(B'/\fm B')\ge \dim B'-\dim R=
\dim A'-\dim R=\dim (A'/\fm A'),
$$
which shows that $x'$ is indeed a generic point of $X_y$.\vspace{3pt}

For the remaining implications we will need:

\begin{slem}\label{flat}
If\/ \textup{\ref{birational}\kern.5pt(ii)}$'$ holds then $f\<\mu$ is flat and the cokernel of\/
$\sO_{\<\<X}\to \mu_*\sO_Z$ is $Y\<\<$-\kern.5pt flat.
\end{slem}

\begin{proof}
Let $x\in X$, $y\set f(x)$, and let  $R\to A\subset B\>$ be as above.
It suffices to show that both $B$ and~$B/A$ are $R$-flat. 
The point $x$ specializes to a closed point of~$X\<$,
whose $f$-image is a specialization of~$y\>$; and for proving flatness we can replace $x$ and $y$ by these specializations. Thus we may assume that $x$ is a closed point, to which (ii)$'$ applies.

Since  $X$ is reduced (by~\ref{factA}(i)), therefore so is $A$. The
assumed birationality of~$\mu$ makes the natural map $A\to B$
injective, and  with $h$ as in Lemma~\ref{hexist}, $B\subset
h^{-1}\<\<A\subset K_A$.

It follows from (ii)$'$ that $h$ is $B_1$-regular;
and since $hB_1\subset A_1$, we deduce that the map 
$\iota_1\colon B_1\to h^{-1}\<\<A_1\cong h^{-1}\<\<A\otimes_R k$,  derived from the inclusion 
\mbox{$\iota\colon B\hookrightarrow h^{-1}\<\<A$} by tensoring over~$R$
with $k\set R/\fm$, is injective.
Since $h^{-1}\<\<A\cong A$ is $R$-flat,  the natural exact sequence\looseness=-1
\begin{equation} \label{heqn}
0=  \tor^R_1(h^{-1}\<\<A,\>k) \to\tor^R_1(h^{-1}\<\<A/\<B,\>k) \to 
B\otimes_R k\xrightarrow{\iota_{\<1}}
 h^{-1}\<\<A\otimes_R k\tag{\ref{hexist}.1}
\end{equation}
gives  $\tor^R_1(h^{-1}\<\<A/\<B,k)=0$, i.e., the finite $A$-module
$h^{-1}\<\<A/B$ is $R$-flat (by the local criterion of flatness \cite[p.\,174,
22.3, ($3'$)$\Rightarrow(1)$]{Ma}); and consequently $B$ is $R$-flat.
Moreover, since (as above) $A_1\subset B_1$,  therefore the natural exact sequence
$$
0=\tor_R^1(B,\>k)\to \tor_R^1(B/A,\>k)\to A_1\to B_1
$$
gives $\tor_R^1(B/A,\>k)=0$, so  $B/A$ is $R$-flat.
\end{proof}

 (ii)$'\Rightarrow{}$(i). There is an easy reduction to
the case where, with $R\to A\to B$ as above, $Y=\spec R$, $X=\spec A$ and $Z=\spec B$ ($f$ and $\mu$ being the obvious maps), and $Y_1=\spec S$ with $S$ a local $R$-algebra. 

\emph{There exists an\/ $A$-regular element\/ $h$ such that $hB\subset A$ and 
$h^{-1}A/\<B$ is\/ $R$-flat.} These assertions were deduced from (ii)$'$ in the proof of~\ref{flat}, so we have them in the case where $x$ is closed in~$X\<$, a case to which, however,  the general case reduces (as in the proof of~\ref{flat}) by specialization.

Now from~\ref{factA}(i) it follows that $A\otimes_R S$ is reduced.
Since, by~\ref{flat}, $B/A$ is $R$-flat, therefore $\tor^R_1(B/A,S)=0$, and
so the natural map $A\otimes_R S\to B\otimes_R S$ is injective. Moreover,
$A/hA\cong h^{-1}\<\<A/A$ is $R$-flat  (take $B=A$ above, i.e., use (ii)$'$ for $Z=X$),  
whence\looseness=-1 
$$
\tor^R_1(A/hA,S)=\tor^R_1(h^{-1}\<\<A/A,S)=0,
$$ 
so that multiplication by $h\otimes 1$ is an injective endomorphism 
of $A\otimes_R S$ (resp.~\mbox{$h^{-1}\<\<A\otimes_R S$}). It follows that 
$h^{-1}\<\<A\otimes_R S$ is isomorphic to an $S$-submodule 
of the quotient ring $K_{A\otimes_R S}$. 
But since $h^{-1}\<\<A/\<B$ is $R$-flat,
the natural map $B\otimes_R S\to h^{-1}\<\<A\otimes_R S$ is injective, and so
$B\otimes_R S$ is isomorphic to a subring of $K_{A\otimes_R S}$. 
This gives (i).\vspace{1pt}

(ii)$'\Rightarrow{}$(iv). By ~\ref{flat}, $f\<\mu$ is flat; and for any $y\in Y\<$, taking  $S$ in the preceding paragraph to be  a finite purely inseparable extension of~$\kappa(y)$ we see that $Z_y$ is geometrically
reduced. So the map $f\<\mu$ is reduced, and hence,  by~\ref{finitefed}, satisfies \cond.
\end{proof}

Here is one situation where the conditions in~\ref{birational} are satisfied.

\begin{sprop}\label{injectivity}
Suppose\/ $Y$ is a $d$-dimensional irreducible regular
scheme and that \mbox{$\triple$} satisfies \cond. If\/ $Z$
as in\/~\textup{\ref{factA}} satisfies the Serre condition\/ $(S_{d+1}),$ $(\>$for example if\/ $Z$ is Cohen-Macaulay$)$ then the
conditions in\/ ~\textup{\ref{birational}} hold. In particular, if\/
$d=1$ and\/ $Z$ is normal then those conditions hold.
\end{sprop}

\begin{proof}
We reduce
as in the proof of~\ref{birational} to considering a
regular local ring $(R,\fm)$ of dimension, say, $n\le d$, a flat
equidimensional local $R$-algebra~$A$, and a finite torsion-free
extension $B\supset A$ satisfying (S$_{n+1}$). 
It is to be shown that any associated prime ideal~$\fP$
of~$\fm B$ intersects $A$ in a minimal associated prime of $\fm A$ (cf.~\ref{birational}(ii)).

Set $\fp\!:=\fP\cap A$. 
Since $B$ is integral over $A$ 
and $B/\fP$ is integral over $A/\fp $,  and since $A$ and~$B$ 
are biequidimensional (see~\ref{finitefed}), 
it holds that
$$
\dim B_\fP\le \dim A_\fp =\dim A-\dim A/\fp =\dim B-\dim B/\fP=\dim B_\fP\>,
$$
so that $(*)\colon \dim B_\fP= \dim A_\fp$. 

Let $(r_1,r_2,\dots,r_n)$ generate $\fm$. Any prime $A$-ideal~$\fq$
has height at least that of its inverse image in $R$ 
\cite[p.\,68, 9.5]{Ma}. Hence for $m\le n$, if
$\fq\supset(r_1,\dots,r_m)A$ 
then \mbox{$\dim A_\fq\ge m$.} Since $B$ satisfies
$(S_{n+1})$,\vspace{.6pt} it follows from $(*)$ and from
\cite[(5.7.5)]{EGA}, by induction on~$m$, that the sequence
$(r_1,\dots,r_m)$ is $B$-regular and every associated prime ideal
$\fQ$ of $(r_1,\dots,r_m)B$ has height $m$.  
Thus $n=\dim B_\fP=\dim A_\fp$, and the conclusion follows.
\end{proof}

%:Thm
 
\begin{thm}\label{Snormal} Let\/ $f\colon X\to Y$ be a reduced
scheme map with\/ $Y$ normal, and $\mu\colon Z\to X$  a finite map. 
If\/ $\mu$ is a
simultaneous normalization of\/ $f$ then $\mu$ is a normalization and 
the nonempty fibers of\/ $f\<\mu$ are geometrically normal.
The converse holds whenever\/ $f$~satisfies\/~$\cond$.
\end{thm}

\begin{proof} Suppose $\mu$ is a simultaneous normalization of
$f$. The fiber $Z_y$ is geometrically normal by definition.
That $X$ is reduced and $Z$ is normal follow from \cite[p.\,184, Cor.\,(ii)]{Ma}; so to show that $\mu$ is a normalization, 
we need only  prove that $\mu$ is birational,  
in other words that $\mu$~induces a bijection from the set of generic points $z\in Z$ to the set of generic points of $X$ such that for each such~$z$, the corresponding local homomorphism $\sO_{\<\<X\<\<,\>\>\mu(z)}\to\sO_{Z\<,\>z}$ is an isomorphism. \vspace{.6pt}
 Flatness of~$f$ and of~$f\<\mu$ implies that every generic point of $X$ and every generic point of $Z$ maps to a generic point of $Y\<$. We may therefore localize at a generic point of $Y\<$, i.e., assume that $Y$ is the spectrum of a field. Then by definition $\mu$ is a normalization map,  thus birational.
\vspace{1pt}
 
For the converse implication,\vspace{-1pt} with $\alpha=\alpha^{}_{Z,\>\kappa(y)}\colon \overline{X_y}\to Z_y\ (y\in f(X))$ as in~\ref{factA}(i),
 2.2.1(iii)$'$ gives that $f\mu$ satisfies \cond, so is flat, and
that
the natural map \smash{$\sO^{}_{\<\<Z_y}\to\alpha_{*}\sO_{\overline{\!X_y\<\!}}\>\>$} is~injective. Normality of~$Z_y$ implies then that this map is surjective too, whence $\alpha$, being finite, is an isomorphism. Thus the map $Z_y\to X_y$ induced by $\mu$ is a normalization (see \ref{factA}(i)). 
\end{proof}

\begin{scor}\label{CorSnormal}
A map\/ $f\colon X\to Y$ satisfying\/ $\cond$ is equinormalizable iff for one\/ $($hence any$)$
normalization map $\nu\colon \BX\to X$  and all\/ $y\in  f(X),$ $\BX_{\<\<y}$ is geometrically normal.
\vspace{1pt}
\end{scor}

\begin{scor}\label{simnorlocal}
Let\/ $f\colon X\to Y$ satisfy $\cond$. Assume that for all\/ $x\in X$ the local ring\/~  $\sO_{\<\<X\<,\>x}$ satisfies\/
$\ff{nor}$. Then\/ $f$~is equinormalizable iff\vspace{.6pt} 
for each\/ $x\in X$
the map\/~$\hat f_{\<x}^{}$ in\/ \textup{~\ref{factA}\kern.5pt (iv)} is equinormalizable. 
\end{scor}

\begin{proof} Let $k$ be a field,\vspace{.4pt} let $(A,\fm)$ be a
reduced local $k$-algebra satisfying $\ff{nor}$,\vspace{.5pt} and
let~$\,\overline{\!A}$ be its integral closure.\vspace{.5pt} By
\cite[(7.6.1)]{EGA}, $\,\overline{\!A}$ is a finite $A$-module,\vspace{-1.7pt} the
completion $\hat A$ is reduced,\vspace{.5pt} and  $\,\overline{\!\hat A}$ is the $\fm$-adic completion of $\,\overline{\!A}$. \vspace{.6pt} 
Moreover, \emph{$A$ is geometrically normal over $k$ if and only if $\hat A$ is.} 
For if the local ring $A'\set A\otimes_k k'$ is normal for every finite purely inseparable 
field extension $k'$ of $k$, then so is its completion $\hat A\otimes_k k'$ 
(since $A'$ satisfies $\ff{nor}$, see~\ref{localgzn}); and conversely,
if $\hat A\otimes_k k'$ is normal then so is $A'=(\hat A\otimes_k k')\cap K_{A'}$. Similar considerations hold if $A$ is semilocal.

Now suppose that $f$ is equinormalizable. Let $C\set\sO_{\<\<X\<,\>x}$ and $D\set\overline C$. As above, $\spec \hat D\to\spec\hat C$ is a normalization map. With $y\set f(x)$ and $\fm$ the maximal ideal of~$\sO_{Y\<,\>y}$,
we have that $\hat D/\fm\hat D$ is the completion of the semilocal
ring $D/\fm D$; and since, by assumption, the latter is geometrically
normal, therefore so is the former. 
Thus  by ~\ref{CorSnormal}, $\hat f_{\<x}^{}$ is equinormalizable.

Conversely, if  $\hat f_{\<x}^{}$ is equinormalizable then $\hat D/\fm\hat D$ is geometrically normal, whence so is~$D/\fm D$. It follows therefore from~\ref{CorSnormal} that if $f$ satisfies \cond\ and $\hat f_{\<x}^{}$ is equinormalizable
 for every $x\in X$ then $f$ is equinormalizable.
\end{proof}

\section{Partially numerical criteria for equinormalization.}

Under suitable conditions, Proposition~\ref{diffdel} 
and Corollary~\ref{regular} give a criterion for equinormalizability
of $f\colon X\to Y$ in terms of constancy of a numerical 
invariant~$\delta$ associated\-
to each fiber~$X_{\<\<y}\set f^{-1}y\ (y\in Y)$, namely, with $\kappa(y)$
the residue field of~$\sO_{\>Y\<\<\<,\>\>y}$, $f_y\colon
X_{\<\<y}\to\spec\kappa(y)$ the obvious map,  and  
$\overline{\<\<\sO_{\<\<X_{\<\<y}}\<\<}\>\>$ the integral closure of~
$\sO_{\<\<X_{\<\<y}}$,
$$
\delta_y\set\dim_{\kappa(y)}f_{y*}(\>\>
\overline{\<\<\sO_{\<\<X_{\<\<y}}\<\<}\>\>/\sO_{\<\<X_{\<\<y}}).
$$
It is assumed here that  $X_{\<\<y}$ has isolated nonnormal points,
with residue
fields finite over~$\kappa(y)$, so that this $\delta$ is finite. 

When $f$ is a flat projective map, we can associate to each fiber
$X_{\<\<y}$ the Hilbert polynomial of
$\>\>\overline{\<\<\sO_{\<\<X_{\<\<y}}\<\<}\>\>/\sO_{\<\<X_{\<\<y}}$; 
and show
under suitable conditions that  equinormalizability is equivalent to
the local constancy of this function of $y$ (see Proposition~\ref{hilbconst}, noting that the Hilbert polynomial of $\sO_{\<\<X_y}$ is locally independent of $y$).%
\footnote{There should be something interesting to be said about this 
criterion being \emph{global} on the
fibers whereas equinormalizability is a \emph{local} condition
(see~\ref{simnorlocal})---but we don't know what that might be.}

It may be noted that when a fiber of a projective map has isolated
nonnormal points, the above Hilbert polynomial is just the constant $\delta$. 

These results will be improved in \S\ref{sec:main}---the above-mentioned
``suitable conditions'' will be weakened to where they refer 
\emph{solely to the fibers themselves.}

\begin{defn} Let $k$ be a field and let $g:X\to \spec k$ be a scheme
map with $X$ a reduced Nagata scheme.
Let $\fC\subset \OX$ be the conductor of 
the normalization \mbox{$\nu\colon\BX\to X\<$,} (a finite map, 
see~\ref{reducednagata}), i.e.,~the annihilator of the
$\OX$-module $\nu_*\sO_{\<\<\BX}/ \OX$; 
and assume that the closed subscheme
$X_\fC\subset X$ corresponding to the coherent $\OX$-ideal~$\fC$
is \emph{finite} over $k$. When these conditions hold 
we say ``$\delta_k(X)$ is finite" and set
$$
\delta_k(X)\set \dim_k g_*(\nu_*\sO_{\BX}/ \OX)
=\sum_{x\in X}\dim_k(\>\>
\overline{\<\<\sO_{\<\<X\!,\>x}\<\<}\>\>/\sO_{\<\<X\!,\>x})<\infty.
$$
If $X$ is affine, 
say $X=\spec A$, we write $\delta_k(A)$ in place of $\delta_k(X)$.
\end{defn}

\begin{defn}\label{spade} 
A ring homomorphism $\phi\colon R\to A$ satisfies \scond\ if\kern.5pt:
 \begin{enumerate}
 \item $(R,\fm,k)$ is a  normal local ring satisfying $\ff{nor}$ and
such that the residue field~$k$ is either of characteristic\/ $0$ 
or of characteristic\/ $>0$ and perfect. 

 \item $A$ is a formally equidimensional Nagata ring.

 \item The map $\phi$ is flat, $\fm A$ is contained in every maximal $A$-ideal,
 $A/\fm A$ is reduced  and $\delta_k(A/\fm A)$ is finite.

\item $\BA/A$ is a finite $R$-module.
 \end{enumerate}
\end{defn}

\pagebreak[3]
\begin{srmks}\label{spadelem}

(a) Suppose  that $R$ is a complete local ring,  or that $R$ is
henselian 
and $A$ is a localization of a finitely generated $R$-algebra, or that  $R$
and~$A$ are both analytic local rings, i.e., homomorphic images of
convergent power series rings over a complete nondiscrete valued  field.
Then conditions (1), (2) and (3) in ~\ref{spade} imply (4).%
\footnote{Complete local rings and analytic local  rings are  
all universally catenary Nagata rings which satisfy~$\ff{nor}$---in fact they are \emph{excellent}
\cite[7.8.3, 5.6.4]{EGA}. For various proofs see \cite[Chap.~0, (22.3.2)]{EGA} plus \cite[p.\,193, (45.5)]{Nag}, \cite[p.\,291, Remark]{M1}, \cite[p.\,96, Satz 3.3.3]{BKKN}, \cite[p.\,1001, Thm.\,2.5]{Kz}, or  the last sentence in the introduction to \cite{Kl}.}

\begin{proof}
 If $\>\fP$ is a prime $A$-ideal  containing $\fm A$ and such that $A_\fP/\fm A_\fP$ is normal, hence geometrically normal over $k$ (since $k$ is perfect or of characteristic~0), then by~\ref{Nishimura} the
homomorphism $R\to A_\fP$ is normal, so $A_\fP$ is normal
\cite[p.\,184, Cor.\,(ii)]{Ma}.\looseness=-1

Finiteness of $\delta_k(A/\fm A)$ means that if $\>\fM$ is a prime $A$-ideal containing $\fm A$ such that $A_\fM/\fm A_\fM$ is \emph{not} normal, then $\fM$ is a maximal ideal and
$[A/\fM:k]<\infty$; and moreover, there are only finitely many such $\fM$.

The ring $A$ is reduced (see~\ref{factA}(i)), so the integral closure $\BA$ is a finite $A$-module (see~\ref{reducednagata}). Let $\fC$ be the 
$A$-conductor, i.e., the annihilator of the $A$-module $\BA/A$. Then
$\fC_\fP=A_\fP$ for any $\fP$ as above, whence, by the preceding paragraph,
$(A/\fC)\otimes_R k$ is a finite-dimensional $k$-vector space. Hence
$A/\fC$ is a finite $R$-module:
\begin{itemize}
\item if $R$ is complete, by \cite[p.\,58, 8.4]{Ma}, since $\fm A$ is contained in the   Jacobson radical of~$A$,   so $A/\fC$ is $\fm$-adically 
separated;

\item if $R$ is henselian and $A$ is a finitely generated $R$-algebra, 
by \cite[18.5.11\;c$')$]{EGA} (since $A/\fC$ has only finitely many 
maximal ideals, all of which contract in $R$ to~$\fm$);\looseness=-1

\item if $R$ and $A$ are analytic local rings, 
 \cite[p.\,18-01, Thm.\,1]{C}.

\end{itemize}
In any of these cases 
the finite $A/\fC$-module $\BA/A$ is also finite over $R$.
\end{proof}

(b) If $\phi\colon (R,\fm,k)\to A$ satisfies \scond\ then
 $\spec\phi\colon\spec A\to \spec R$ satisfies ~\cond,
 and the \/ $\fm$-adic completion $\>\hat\phi\colon \hat R\to\hat A\>$ 
satisfies\/ \scond. If, in addition, $A$ satisfies $\ff{nor}$ and 
$A/\fm A$~has finite Krull dimension, then $\hat A$ satisfies~$\ff{nor}$.

\begin{proof}  Since $k$ is perfect or of characteristic 0, 
the $k$-algebra $A/\fm A$ is geometrically reduced
(see~\ref{geornr}).  If $\>\fM$ is any maximal
$A$-ideal then the composition $R\to A\to A_\fM$ is reduced 
(see~\ref{Nishimura}). It
follows that $\spec\phi$ is reduced (since being reduced is a local
property, and every prime $A$-ideal is contained in some~$\fM$). It is now immediate that  $\spec\phi$ satisfies ~\cond.

Since $R$ is normal and satisfies $\ff{nor}$, therefore 
$\hat R$ is normal \cite[p.\,184, Cor.\,(ii)]{Ma}, and of course $\hat R$ satisfies $\ff{nor}$. That $\spec\hat\phi$ satisfies \cond\ is given
by~\ref{factA}(iv) with $I=\fm$ and
$J=\fm A$. Since $\hat A/\fm \hat A=A/\fm A$ and $\hat R$ is complete, 
it follows from (a) that $\hat\phi$ satisfies~\scond.\looseness=-1

In particular, $\hat A$ is universally catenary; and
consequently, since $\fm\hat A$ is contained in every maximal ideal 
of $\hat A$ and $\hat A/\fm \hat A=A/\fm A$ has finite Krull
dimension, therefore  $\hat A$ has finite Krull dimension.
Then Satz 2 of \cite{BR} shows that $\hat A$ satisfies~$\ff{nor}$. 
\end{proof}

\end{srmks}

\begin{prop} \label{diffdel}
Let\/ $\phi\colon (R,\fm,k)\to A$ satisfy \scond. Set\/ $K\set K_R\>,$ 
$A_0\set A\otimes_R K\<,$
$A_1\set A/\fm A,$ $B\set\BA,$  $B_0\set B\otimes_R K\<,$
$B_1\set B/\fm B$.   For\/ $\fp\in\spec R$ 
set\/ $A_{\langle\fp\rangle}\set A\otimes_R \kappa(\fp)$  and\/ 
$B_{\langle\fp\rangle}\set B\otimes_R\kappa(\fp)$.  Then$\>:$\vspace{1pt}
\vspace{2pt}

{\rm(i)} If\/ $f\set\spec\phi$ is equinormalizable then\/ 
$\delta_{\kappa(\fp)}(A_{\langle\fp\rangle})$ is independent of\/~$\fp$.\vspace{1pt}

 Assume further that\/ the map\/ \mbox{$\alpha\colon B_1\to\overline{A_1}$} arising from~\textup{\ref{factA}(i)} is injective. Then$\>:$\vspace{1pt}

{\rm(ii)} If\/ $\delta_{\kappa(\fp)}(A_{\langle\fp\rangle})<\infty$ then
$$
\delta_{\kappa(\fp)}(A_{\langle\fp\rangle})-\delta_{\kappa(\fp)}(B_{\langle\fp\rangle})=\delta_K(A_0).\\[1pt]
$$

{\rm(iii)} If\/ $\delta_k(A_1)\le\delta_K(A_0)$ then
$f$ is equinormalizable.
\end{prop}

\begin{proof} For any simultaneous
normalization $\mu\colon Z\to X\set\spec A$ of $f$, $Z\cong \spec B$  and condition~\ref{birational}(ii) is satisfied 
(see Theorem~\ref{Snormal}), whence so is ~\ref{birational}(i), so that
$\alpha\colon B_1\to\overline{A_1}$ is injective. Thus we may assume this injectivity throughout the proof. This assumption implies condition (ii)$'$ in Proposition~\ref{birational}, so by that Proposition both $B$ and 
$B/A$ are flat $R$-modules; and being finitely generated (by $\scond$), 
\emph{$B/A$ is a finite-rank free\/ $R$-module.} 

We begin with the case $\fp=\fm$. 
By assumption, $B_1$ and its subring $A_1$ have the same integral
closure, 
whence
$$
\delta_k(A_1)-\delta_k(B_1)=\dim_k(B_1/A_1)=\dim_k (B/A\otimes_R k).
$$
Also, it is easy to see that $B\otimes_R K= \overline{A\otimes_R K}$. Therefore,
$$
\dim_k (B/A \otimes_R k)
=\dim_K (B/A\otimes_R K)
=\delta_K (A_0).
$$
Thus (ii) holds for $\fp=\fm$. 
\vspace{1pt}

%:3.3

For arbitrary $\fp$, the implication (ii)$'\Rightarrow{}$(i)  in
Proposition~\ref{birational} shows that $B_{\langle\fp\rangle}$ and
its subring $A_{\langle\fp\rangle}$ have the same  integral closure.  
So one can localize at $\fp$ and argue as before to prove (ii).
Moreover,  if $f$ is equinormalizable then $B_{\langle\fp\rangle}$ is normal, and hence 
$$
\delta_{\kappa(\fp)}(A_{\langle\fp\rangle})=\dim_{\kappa(\fp)}(B_{\langle\fp\rangle}/A_{\langle\fp\rangle})=\textup{rank}_R(B/A)
$$
is independent of $\fp$, proving (i). Finally, if
$\delta_k(A_1)\le\delta_K(A_0)$ then by (ii),
$\delta_k(B_1)=0$, i.e., $B_1$ is normal, hence geometrically normal,
since $k$ is perfect or of characteristic 0; and since any prime $B$-ideal
is contained in a maximal ideal, which contracts to a maximal ideal in~$A$, hence to $\fm$ in $R$,
it follows from~\ref{Nishimura} that all the fibers of 
$\spec B\to\spec R$ are geometrically normal.  
Thus Corollary~\ref{CorSnormal} gives (iii).
\end{proof}

\begin{scor}\label{regular}
Let\/  $\phi\colon R\to A$ satisfy\/ $\scond,$ and assume further that $R$ is a regular local ring of dimension, say, $d$ and that $\BA$ satisfies the Serre condition\/ $(\textup{S}_{d+1})$ $($which it always does in case\/ $d=1).$ Then with the notation of Proposition\/~\textup{\ref{diffdel},} $\spec\phi$ is equinormalizable iff\/ $\delta_k(A_1)=\delta_K(A_0)$.

\begin{proof}
In view of Proposition~\ref{injectivity}, this results from Proposition~\ref{diffdel}.
\end{proof} 
\end{scor}
\pagebreak[3]

%:3.4

We turn now to the case where the map $f\colon X\to Y$, satisfying
\cond, is \emph{projective}. 

We need some notation.  For any $Y$-scheme $W$ and any
$y\in Y\<$, with residue field $\kappa(y)$, let  
$W_{\<y}\set W\otimes_Y\spec\kappa(y)$  be the fiber over $y$, and let 
$\iota^W_y\colon W_{\<y}\to W$ be the projection.
For any $\sO_W$-module
$\cF\>$ let $\cF_y$ be the $\sO_{W_{\<y}}$-module $(\iota^W_y)^*\cF$.
\vspace{.6pt}

Projectivity of $f\colon X\to Y$ entails the existence of invertible
$\sO_{\<\<X}$-modules
which are very ample relative to $f$. Fix one such and call it $\cL$.
Then $\cL_y$ is very ample relative to the projection
$f_y\colon X_{\<y}\to \spec\kappa(y)$ \cite[(4.4.10)(iii)]{EGII}.
 For any coherent $\sO_{\<\<X_{\<\<y}}$-module\/~ 
$\mathcal M$ set $\mathcal M(n)\set \mathcal M\otimes_{\sO_{\<\<X_{\<\<y}}}\cL_y^{\otimes n}$ and
let $\sH_y(\mathcal M)$ be the polynomial function (depending on~
$\cL$) which
takes integers $n\gg 0$ to 
$\dim_{\kappa(y)}\tH^0(X_{\<\<y},\>\mathcal M(n))$.   
(See \cite[(2.5.3)]{EGIII}.)

\begin{prop}\label{hilbconst}
Let\/ $(R,\fm,k)$ be a  normal local ring satisfying $\ff{nor}$ and
such that\/ $k$ is either of characteristic\/ $0$ 
or of characteristic\/ $>0$ and perfect.  Let\/ \mbox{$f\colon X\to
Y\set\spec R$} be a projective map whose fibers are all geometrically reduced and 
which is locally equidimensional \textup{\cite[Err$_{\textup{IV}}$, 35]{EGA}}.%
\footnote{Pub.\ Math.\ I.H.E.S.\ {\bf 32}, p.\,357.}
 Let $\cL,$ $\sH_y$ be as above. Let $\mu\colon Z\to X$ be a normalization map. 
Let\/ $y_0$ and\/ $y_1$
be, respectively, the closed and generic points of\/~$Y\<,$ and set\/
$X_i\set X_{\<y_i},$ resp.~$Z_i\set Z_{y_i}.$ Then$\>:$\vspace{1pt}

\textup{(i)} If   
the map\/ \mbox{$\alpha\colon\overline{X_{\<1}}\to Z_1$} 
from~\textup{\ref{factA}(i)} is schematically dominant \textup{(\S\ref{dominant}),} and if
$$
\sH_{y^{}_1}\<\<(\>\>\overline{\<\<\sO_{\<\<X_{\<1}}\<\<}\>\>)=
\sH_{y^{}_0}\<\<(\>\>\overline{\<\<\sO_{\<\<X_0}\<\<}\>\>)
$$ 
then\/ $f$ is equinormalizable.\vspace{1pt}

\textup{(ii)} If\/ $f$ is equinormalizable then\/ 
$\sH_y(\>\>\overline{\<\<\sO_{\<\<X_{\<\<y}}\<\<}\>\>)$ 
is independent of\/~$y$.
  \end{prop}

\begin{proof} First of all, $f$ satisfies \cond. In fact:

\begin{slem}\label{conds} 
Let\/ $(R,\fm,k)$ be a  normal local ring satisfying $\ff{nor}$. Let 
$f\colon X\to \spec R$ be a finite-type scheme map whose fibers are all geometrically
reduced. Then the following conditions are equivalent.

{\rm (i)} $f$ satisfies \cond.

{\rm (ii)} $f$ is locally equidimensional.

{\rm (iii)} For each\/ $x\in X,$  $\sO_{\<\<X\<,\>x}$ is equidimensional, and\/ 
$f$ is universally open.

{\rm (iii)$'$} For each\/ $x\in X,$ $\sO_{\<\<X\<,\>x}$ is equidimensional, and\/ 
$f$ is open.

{\rm (iv)} For each\/ $x\in X,$ $\sO_{\<\<X\<,\>x}$ is
equidimensional, and\/ $f$ is flat.

{\rm (v)} For each\/ $x\in X$ and $y=f(x)$, $\sO_{\<\<X\<,\>x}$ is equidimensional,
and with\/ $\fm_y$ the maximal ideal of\/~$\sO_{Y\!,\>\>y},$ 
$$
\dim \sO_{\<\<X\<,\>x}=\dim\sO_{Y\!,\>\>y}+\dim(\sO_{\<\<X\<,\>x}/\fm_y\sO_{\<\<X\<,\>x}).
$$
\end{slem}

\begin{proof} (i)${}\Rightarrow{}$(iv). Trivial.

(iv)${}\Rightarrow{}$(i). It suffices to show that $\spec R$ is a universally catenary
Nagata scheme, since then the same will be true of $X$ (see~\ref{unicat}, \ref{localnagata}). Since $\hat R$ is normal \cite[p.\,184, Cor.\,(ii)]{Ma}, therefore $R$ is universally catenary (see~\ref{qunmixed}, \ref{lqunmixed}). Since $R$ satisfies $\ff{nor}$,
therefore $R$ is a Nagata ring (see~\ref{univjap}).

(iv)${}\Rightarrow{}$(iii). \cite[(2.4.6)]{EGA}.

(iii)${}\Rightarrow{}$(iv). \cite[(15.2.3)]{EGA}.

(iii)${}\Leftrightarrow{}$(iii)$'$.  \cite[(14.4.3)]{EGA}. (Normal schemes are
geometrically unibranch \cite[(6.15.1)]{EGA}.)

(iv)${}\Rightarrow{}$(v).  \cite[(6.1.2)]{EGA}.

(v)${}\Rightarrow{}$(ii). \cite[(13.3.6)]{EGA}.

(ii)${}\Rightarrow{}$(iii). \cite[(14.4.4)]{EGA}.
\end{proof}

Now, as in the proof of Proposition~\ref{diffdel}, 
we may assume throughout that
$\alpha$ is schematically dominant.\vspace{1pt}

The projective map $f$ takes closed points of $X$ to the closed point of $Y\<$.
The assumed dominance of $\alpha$ implies, via~\ref{birational}, that $f\<\mu$ is flat, and so $\overline{\sO_{\<\<X}}=\mu_*\sO_Z$ is $Y\<$-flat. Hence by \cite[7.9.13]{EGIII}, $\sH_y(\<\<(\>\>\overline{\<\<\sO_{\<\<X}\<\<}\>\>)_y)$ is independent of $y$. 

If $f$ is equinormalizable, then $\BX_y=\overline{X_y}$, and one deduces (directly, or by \cite[p.\,366, (9.3.3)]{EGA1}) a natural isomorphism
$\<\<(\>\>\overline{\<\<\sO_{\<\<X}\<\<}\>\>)_y\iso\>\>\overline{\<\<\sO_{\<\<X_{\<\<y}}\<\<}\>\>$, proving (ii). 

As for (i), since $\alpha$ is schematically
dominant, therefore
$ \overline{\<\<\sO_{\<\<X_{\<1}}\<\<}\>\>
\supset(\>\overline{\<\sO_{\<\<X}\!}\>\>)_{\<y^{}_1}$, and so
\begin{align*}
\sH_{y^{}_1}\<\<(\>\>\overline{\<\<\sO_{\<\<X_{\<1}}\<\<}\>\>)
  &=
  \sH_{y^{}_1}\<\<\Big(\>\>\overline{\<\<\sO_{\<\<X_{\<1}}\<\<}\>\>
  /(\>\overline{\<\sO_{\<\<X}\!}\>\>)_{\<y^{}_1}\Big) +
\sH_{y^{}_1}\<\<\Big((\>\overline{\<\sO_{\<\<X}\!}\>\>)_{\<y^{}_1}\Big) \\
 & =
  \sH_{y^{}_1}\<\<\Big(\>\>\overline{\<\<\sO_{\<\<X_{\<1}}\<\<}\>\>
  /(\>\overline{\<\sO_{\<\<X}\!}\>\>)_{\<y^{}_1}\Big) +
\sH_{y^{}_0}\<\<\Big((\>\overline{\<\sO_{\<\<X}\!}\>\>)_{\<y^{}_0}\Big) \\
   & =
  \sH_{y^{}_1}\<\<\Big(\>\>\overline{\<\<\sO_{\<\<X_{\<1}}\<\<}\>\>
  /(\>\overline{\<\sO_{\<\<X}\!}\>\>)_{y^{}_1}\Big) +
\sH_{y^{}_0}\<\<\Big(\>\>\overline{\<\<\sO_{\<\<X_0}\<\<}\>\>\Big). 
\end{align*}
The equality $\sH_{y^{}_1}\<\<(\>\>\overline{\<\<\sO_{\<\<X_{\<1}}\<\<}\>\>)=
\sH_{y^{}_0}\<\<(\>\>\overline{\<\<\sO_{\<\<X_0}\<\<}\>\>)$ implies then that 
$\sH_{y^{}_1}\<\<(\>\>\overline{\<\<\sO_{\<\<X_{\<1}}\<\<}\>\>
  /(\>\overline{\<\sO_{\<\<X}\!}\>\>)_{\<y^{}_1})
  =0$, whence $\>\>\overline{\<\<\sO_{\<\<X_{\<1}}\<\<}\>\>=
(\>\overline{\<\sO_{\<\<X}\!}\>\>)_{\<y^{}_1} $, i.e., 
$(\>\overline{\<\sO_{\<\<X}\!}\>\>)_{\<y^{}_1}$ \emph{is normal.}  
(Indeed, for any coherent $\sO_{\<\<X_{\<\<y}}$-module\/ 
$\mathcal M$, $\sH_y(\mathcal M)=0 \Leftrightarrow \mathcal M=0$ 
because for all $n\gg0$, $\mathcal M(n)$ is generated by global
sections.) 
That $f$ is equinormalizable follows now just as in the 
proof of~\ref{diffdel}(iii).
\end{proof}

As before, Proposition~\ref{injectivity} yields:
\begin{scor}\label{regular2}
In\/~\textup{\ref{hilbconst}} assume  that $R$ is a regular local ring of dimension, say, $d$ and that $Z$ satisfies the Serre condition\/ $(\textup{S}_{d+1})$ $($which it always does in case\/ $d=1).$ Then  $\alpha$ has to be  schematically dominant, and so $f$ is equinormalizable iff\/
$\sH_{y^{}_1}\<\<(\>\>\overline{\<\<\sO_{\<\<X_{\<1}}\<\<}\>\>)=
\sH_{y^{}_0}\<\<(\>\>\overline{\<\<\sO_{\<\<X_0}\<\<}\>\>).$ 
\end{scor}

\section{Main Theorems.}\label{sec:main}

Here is the first main theorem, affirming that for families of curves  it is not necessary in Proposition~\ref{diffdel}(iii) to assume that the map $\alpha$ is injective.

\begin{thm}\label{main}
 Let\/ $\phi\colon (R,\fm,k)\to A$ satisfy\/ $\scond.$ Suppose also that\/ $R$ and\/ $A$ satisfy one of the conditions in Remark\/~\textup{\ref{spadelem}(a)} 
and that\/ $A/\fm A$ has Krull dimension\/~$1.$  Set\/ $K\set K_R\>,$ 
$A_0\set A\otimes_R K\<,$ $A_1\set A\otimes_R k.$ If\/ $\delta_k(A_1)=\delta_K(A_0)$
then\/ $f\set\spec\phi$ is equinormalizable.
\end{thm}

\begin{proof}
Set $B\set\BA,$  $B_0\set B\otimes_R K\<,$
$B_1\set B\otimes_R k$. 
By Proposition~\ref{diffdel}, it suffices to show that\vspace{-1pt}
 the map
\mbox{$\alpha\colon B_1\to\>\>\overline{\<\<A_1\<\<}\>\>$} given by~\textup{\ref{factA}(i)} is injective.
\vspace{2pt}

Recall that $f$ satisfies \cond\ (Remark~\ref{spadelem}(b)), so that by~\ref{factA} both~$A$ 
and~$A/\fm  A$ are reduced, and for any minimal prime divisor $\fP$ of~$\fm A$ the local ring $A_\fP$ is normal. This being so, the proof of Lemma~\textup{\ref{hexist}} is valid in the present situation.

\begin{slem}\label{torsionfree}
If\/ $h$ is as in Lemma~\textup{\ref{hexist}} then the $R$-module\/ $h^{-1}\<A/B$ is torsion-free.  
%and it is flat if and only if\/ $\alpha$ is injective.
\end{slem}

\begin{proof} Multiplication by the $A$-regular element $h$ is an $R$-isomorphism
$h^{-1}\<A/B\cong A/hB$. Note that the $A$-ideal $hB$ is the integral closure of $hA$.
As $A$ is formally equidimensional,
\cite[p.\,189, Thm.\, 2.12]{R3}
implies that every associated prime $A$-ideal $\fQ$ of $hB$ is minimal.

Now these $\fQ$ are also the minimal prime
 divisors of $hA$, because $hA \subset hB \subset \sqrt{hA}$. Since $A$ is $R$-flat and multiplication by $h$ is an injective endomorphism of 
 $A/\fm A=A\otimes_R k$, therefore $\tor^R_1(A/hA,k)=0$; so by \cite[p.\,174, 22.3]{Ma}, $A/hA$  is $R$-flat, hence torsion-free, and therefore $\fP \cap R =0 $ for every associated prime $\fP$ of $hA$, in particular  for $\fP=\fQ$. Hence $A/hB$ is $R$-torsion-free. 
\end{proof}

\begin{srmk} \label{torfree}
As just shown, $h^{-1}A/A\cong A/hA$ is $R$-flat; and since 
$A\subset\BA\subset h^{-1}A$, therefore $\BA/A$ \emph{is\/ $R$-torsion-free.}
\end{srmk}

Note next that the above $h$ can be chosen so that $A/hA$ \emph{is a
finite $R$-module.}  Indeed, as in Remark~\ref{spadelem}(a), it
suffices to arrange that $A/hA\otimes_R k$ be finite over $k$. Since
$A/\fm A$ has dimension~one and $h$ is $A/\fm A$-regular, therefore
$A/hA\otimes_R k$ is artinian, so it suffices that every maximal
$A$-ideal~$\fM$ containing $h$ satisfy $[A/\fM:k]<\infty$; and again as in Remark~\ref{spadelem}(a), this will be so if $A_\fM/\fm
A_\fM$ is not normal. Reviewing the proof of Lemma~\ref{hexist}, one
finds it enough to show that \emph{if\/ $A_\fM/\fm A_\fM$ is normal
then neither\/~$\fC'$ nor\/ $\fC$ is contained in\/ $\fM$.} For $\fC'$ this is evident; 
and for $\fC$
it follows from~\ref{Nishimura}---which with \cite[p.\,184, Cor.\,(ii)]{Ma} shows that $A_\fM$ is normal, so that $\fC A_\fM\>$, the conductor of~$A_\fM\>$, is the unit ideal.

As observed in the proof of~\ref{torsionfree}, $h^{-1}\<A/A\cong A/hA$ is $R$-flat, 
and hence, being finitely generated, it is $R$-free, of rank, say, $d$. 
For any $R$-module~ $C$ and any $R$-algebra $T$, set 
$C_{\>T}\set C\otimes_R T$.  The natural exact sequence
\addtocounter{equation}{2}
\begin{equation}\label{htor}
0=\tor_R^1(A/hA, T)\to A_{\>T}\xrightarrow{h} A_{\>T}
\end{equation}
shows that $h$ is $A_{\>T}$-regular, so that  there are natural isomorphisms
$$
(h^{-1}\<A)_{\>T}\underset{h\otimes 1}\iso A_{\<T}
\underset{h^{-1}}\iso h^{-1}\<A_{\>T}\subset  K_{A_{\>T}}.
$$

Assume henceforth that $R\subset T\subset K$ and that $T$ is normal and essentially of finite type over~$R$. Then with $\phi_T\colon T\to A_T$ the map induced by $\phi$, ~\ref{factA}(iii) shows that $\spec \phi_T$ satisfies \cond. Hence, by Remark~\ref{torfree}, $\>\>\overline{\<\<A_T\<\<}\>\>/A_T$ is $T$-torsion free. Consequently, 
\begin{equation}\label{hBA}
A_T\subset\>\>\overline{\<\<A_T\<\<}\>\>\subset h^{-1}A_T\>;
\end{equation}
indeed, $A_K=A_T\otimes_T K$ and $\>\>\overline{\<\<A_K\<\<}\>\>=\>\>\overline{\<\<A_T\<\<}\>\>\otimes_T K$, and
since $h\BA\subset A$\vspace{.6pt} therefore $h\>\>\overline{\<\<A_K\<\<}\>\>\subset A_K$, so
for any $a\in\>\>\overline{\<\<A_T\<\<}\>\>$ there is a nonzero $t\in T$ such that
$tha\in A_T$, whence, by torsion-freeness of  $\>\>\overline{\<\<A_T\<\<}\>\>/A_T$,
$ha\in A_T$.
We denote the inclusion $\>\>\overline{\<\<A_T\<\<}\>\>\hookrightarrow h^{-1}A_T$
by $\iota_T^{}$.

With $\delta\set\delta_K(A_1)=\delta_k(A_0)$, 
let $g\colon\bG\to\spec R$ be the Grassmannian of 
locally free rank-$(d-\delta)$ quotients of $h^{-1}\<A/A$ (\cite[p.\,384, (9.7.5)]{EGA1}).  There is then a functorial bijection between $R$-morphisms
$\spec T\to\bG$ and $T$-submodules $L\subset h^{-1}\<A_T/A_T$ such that 
the $T$-module $(h^{-1}\<A_T/A_T)/L$ is locally free of rank $d-\delta$. The map~$g$ is projective (\cite[p.\,390, (9.8.4)]{EGA1}).

Let $\psi\colon\spec K\to\bG$  be the $R$-morphism corresponding to the $\delta$-dimensional $K$-vector space
$B_0/A_0\subset h^{-1}\<A_0/A_0$. 
Regarding $\psi$ as a rational $R$-morphism $\psi_T^{}\colon \spec T\to \bG$ (see \cite[p.\,345, (8.1.11)]{EGA1}), \emph{suppose that the domain of definition\/ $\textup{D}(\psi_T^{})$\vspace{.4pt}
 is all of} $\spec T$.
Then $\psi_T^{}$ corresponds to a $T$-module 
$L_1\subset E\set  h^{-1}\<A_{\>T}/A_{\>T}$ with $E/L_1$ \vspace{.6pt}
locally free of rank $(d-\delta)$, such that $L_1\otimes_T K =B_0/A_0=L_2\otimes_T K$ where $L_2\set \overline{\<\<A_T}/A_T$.\vspace{.6pt} Since $\spec A_T\to\spec T$ satisfies~\cond\ (see Remark~\ref{spadelem}(b) and~\ref{factA}(iii)),\vspace{.6pt} the $T$-module $E/L_2$ is torsion-free (see ~\ref{torsionfree}).
The following Lemma~\ref{uniquetorsion} 
shows then that $L_1=L_2$, so that\vspace{.4pt}
$h^{-1}\<A_T/\>\>\overline{\<\<A_T\<}\>\cong E/L_2=E/L_1$ is
locally free of rank $d-\delta$, and there is a 
 split-exact sequence of $T$-modules
\begin{equation}\label{split}
0\to \>\>\overline{\<\<A_T\<}\>\xrightarrow{\iota_T^{}} h^{-1}\<A_{\>T}\to h^{-1}\<A_{\>T}/\>\>\overline{\<\<A_T\<}\>\to 0.
\end{equation}

%:ici

\addtocounter{slem}{3}
 \begin{slem}\label{uniquetorsion}
Let\/ $ S$ be a commutative domain,  $E$ a torsion-free\/ $S$-module, 
and\/ $L_1,$ $L_2$ $S$-submodules of\/ $E$. If\/  $E/ L_1$
and\/ $E/ L_2$ are both\/ torsion-free and if the natural images of\/
$ L_1\otimes_{ S}K_{ S}$ and\/ $ L_2\otimes_{ S}K_{ S}$ in\/
$ E\otimes_{ S}K_{ S}$ coincide then\/ $ L_1= L_2$.
 \end{slem}

 \begin{proof}
    The natural map $\rho\colon E \to E\otimes_{ S}K_{S}$ is injective. 
    Since the images of $ L_1\otimes_{ S}K_{ S}$
   and $ L_2\otimes_{ S}K_{ S}$ in $ E\otimes_{ S}K_{ S}$ 
   coincide, there exists for each $f\in L_1$ a nonzero $s\in S$ and a $g\in L_2$ such
   that $\ds \rho(f)={\rho(g)}/{s}$, i.e., $\rho(s f)=\rho(g)$, i.e., $s f=g$. Since
   $ E/ L_2$ is torsion-free, $s f=g$ implies that
    $f \in  L_2$, and thus $ L_1 \subset  L_2$. Similarly, 
   $ L_2 \subset  L_1$. 
 \end{proof}

Now if $T=R$, so that $A_T=A$ and $\>\>\overline{\<\<A_T\<}\>=B$, then applying $\otimes_R k$ to the split-exact sequence~\eqref{split} we get an exact sequence
$$
0\to B_1\xrightarrow{\alpha'} h^{-1}\<A_1\to h^{-1}\<A_1/B_1\to 0.
$$
As $\>\>\overline{\<\<A_1\<}\>\subset h^{-1}\<A_1\subset K_{A_1}$, one sees that 
$\alpha'(B_1)\subset \>\>\overline{\<\<A_1\<}\>$,\vspace{.5pt} and that the natural composition
\smash{$A_1\to B_1\xrightarrow[\under{3.75}{\alpha'}]{} \>\>\overline{\<\<A_1\<}\>$} is a normalization map,
whence, by~\ref{factA}(i), $\alpha'=\iota_k^{}\smcirc\alpha$. Thus injectivity of $\alpha'$ implies that of $\alpha$.\vspace{1pt}

In summary:
\emph{For\/ $\alpha$ to be injective 
it suffices that\/ $\textup{D}(\psi_R^{})$ be all of}
$\>\spec R$,  which we will now show to be the case.

\smallskip
Let $Z$ be the schematic closure
of $\psi(\spec K),$ so that $g$ induces a birational projective map
$\gamma\colon Z\to \spec R$.
According to \cite[p.\,347, (8.2.7)]{EGA1}, $\psi$ is defined on 
all of $\>\spec R$ if
$\gamma$ is an isomorphism, for which, since $R$ is normal, it suffices by Zariski's Main
Theorem (see e.g., \cite[(4.4.8)]{EGIII}) that  the 
closed fiber $\gamma^{-1}\{\fm\}$ be zero-dimensional. We need only show then that \emph{$\gamma^{-1}\{\fm\}$ has a unique closed point.}

Let $z$ be any closed point in $Z$. Let $S$ be the local ring of the generic point of
a component of the closed fiber of the normalization of the blowup of the maximal ideal of
$\sO_{Z,z}$. Then $S$ is a discrete valuation ring with fraction
field~$K\<$, essentially of finite type over~$R$ (because $R$ is a
Nagata ring, see~\ref{univjap}), and with
residue field $k_S$  a separable extension of~$k$. 

As above (with $T=S$), the natural map $\spec S\to\bG$ corresponds to $L_2\set\>\>\overline{\<\<A_S\<}\>/A_S$, which, being the kernel of the surjective map 
of free $S$-modules $h^{-1}\<A_S/A_S\to h^{-1}\<A_S/\>\>\overline{\<\<A_S\<}\>$ of respective ranks $d$ and $d-\delta$,  is  free of rank $\delta$. So  
$$
\dim_{k_S}\Bigl((\>\>\overline{\<\<A_S\<}\>\otimes_S k_S)/A_{k_S}\Bigr)=
\dim_{k_S}\Bigl((\>\>\overline{\<\<A_S\<}/A_S)\>\otimes_S k_S\Bigr)
=\delta.
$$ 
As $k_S$ is separable over $k$, \mbox{$\>\>\overline{\<\<A_{k_S}\<}\>=\>\>\overline{\<\<A_1\<}\>\otimes_k k_S$} \cite[(6.14.2)]{EGA}. So
$$
\dim_{k_S}\Bigl(\>\>\overline{\<\<A_{k_S}\<}\>/A_{k_S}\Bigr)=
\dim_{k_S} \Bigl((\>\>\overline{\<\<A_1\<}\>/A_1)\otimes_k k_S\Bigr)
=\delta_k(A_1)=\delta.
$$
Moreover, as above, $\spec A_S\to\spec S$ satisfies \cond,
so \ref{injectivity} gives that the natural map 
$\>\>\overline{\<\<A_S\<}\>\otimes_S k_S\to\>\>\overline{\<\<A_{k_S}\<}\>$
is injective, whence so is the resulting map
$(\>\>\overline{\<\<A_S\<}\>\otimes_S k_S)/A_{k_S}\to\>\>\overline{\<\<A_{k_S}\<}\>/A_{k_S}$. Since the source and target of this last map have the same dimension $\delta$, it is \emph{bijective,} and so there are natural identifications
$$
(\>\>\overline{\<\<A_S\<}\>/A_S)\otimes_S k_S=
(\>\>\overline{\<\<A_S\<}\>\otimes_S k_S)/A_{k_S}=
\>\>\overline{\<\<A_{k_S}\<}\>/A_{k_S}=
(\>\>\overline{\<\<A_1\<}\>/A_1)\otimes_k k_S
\subset h^{-1}\<A_{k_S}/A_{k_S}.
$$
This means that if $\psi_1\colon k\to\bG$ corresponds to 
$\>\>\overline{\<\<A_1\<}\>/A_1\subset h^{-1}\<A_1/A_1$ then the following 
natural diagram commutes:
$$
\begin{CD}
\spec k_S @>>>\spec S \\
@VVV @VVV \\
\spec k @>>\psi_1>\bG
\end{CD}
$$
Hence  $z=\psi_1(\spec k)$, and so $\gamma^{-1}\{\fm\}$ does indeed have a unique closed point.
\end{proof}

\smallskip 
%:(b)'

\begin{thm}\label{main2} Let\/ $(R,\fm,k)$ be a  normal local domain satisfying $\ff{nor}$ and
such that\/ $k$ is either of characteristic\/ $0$ 
or of characteristic\/ $>0$ and perfect.   
Let\/ \mbox{$f\colon X\to Y\set\spec R$} be a projective map whose fibers are all geometrically reduced  and 
which is locally equidimensional. 
Let\/ $\mathcal L$ be an invertible\/ $\OX$-module which is very ample for\/~$f,$ 
and for\/ $y\in Y$ let\/ $\sH_y$ denote the corresponding Hilbert polynomial on the fiber\/ $X_y$ $($see the paragraphs preceding Proposition~\textup{{\ref{hilbconst}}}$)$. 
Let\/ $y_1$ and\/ $y_0$
be, respectively, the closed and generic points of\/~$Y\<\<,$  and set\/
$X_i\set X_{\<y_i}.$  In this situation, one has$\>:$\vspace{2pt}

\textup{(i)} If\/ 
$
\sH_{y^{}_1}\<\<(\>\>\overline{\<\<\sO_{\<\<X_1}\<\<}\>\>)=
\sH_{y^{}_0}\<\<(\>\>\overline{\<\<\sO_{\<\<X_0}\<\<}\>\>)
$
then\/ $f$ is equinormalizable.\vspace{1pt}

\textup{(ii)} If\/ $f$ is equinormalizable then\/ 
$\sH_y(\>\>\overline{\<\<\sO_{\<\<X_{\<\<y}}\<\<}\>\>)$ 
is independent of\/~$y$.

\end{thm}

\begin{proof}
\noindent The proof is analogous to that of
Theorem~\ref{main}. Recall that $f$ satisfies $\cond$, (see
 ~\ref{conds}). In view of Proposition~\ref{hilbconst}, we need only prove (i), for which it suffices to
show, with $Z\to X$ a normalization map and $Z_1\set Z_{y_1}$, that
\emph{if\/
$\sH_{y^{}_1}\<\<(\>\>\overline{\<\<\sO_{\<\<X_1}\<\<}\>\>)=
\sH_{y^{}_0}\<\<(\>\>\overline{\<\<\sO_{\<\<X_0}\<\<}\>\>)$ 
then the map\/ \mbox{$\alpha\colon\overline{X_1} \to Z_1$}
from~\textup{\ref{factA}(i)} is schematically dominant.}\vspace{1pt} 

First, some notation. Let $T$ be an $R$-algebra whose only idempotents are 0 and 1 (i.e., $\spec T$ is connected).
Set $X_T\set X\otimes_R T$.  Hilbert
polynomials on the fibers of the projection 
$f_T\colon X_T\to\spec T$ are defined
via the very ample (relative to $f_T$) invertible
$\sO_{\<\<X_T}$-module $\cL_T\set \cL\otimes_R T$ (cf.~again,  the paragraphs preceding Proposition~\textup{{\ref{hilbconst}}}).
If $\cF$ is a $T$-flat coherent $\sO_{\<\<X_T}$-module then by
\cite[(7.9.13)]{EGIII} the Hilbert polynomial
$\sH_t(\cF_t)$ is the same for all $t\in\spec T$. We denote that polynomial
by $\sH_T(\cF)$.\vspace{1pt}

We will need some global object to take the place of the element $h$ in Lemma~\ref{hexist}. By \cite[(4.4.7)]{EGII} we may assume there exists a graded $R$-algebra $A=R\oplus A_{[1]}\oplus A_{[2]}+\cdots$ generated by  a finite $R$-module $A_{[1]}$, such that $X= \proj A$ and $\cL=\OX(1)$. We may also assume that we are not in the trivial situation 
where $X_1$---and hence $X$---is empty. 
Let $\mathcal C$ be the conductor of $\OX$, and let $\fC$ be a graded
$A$-ideal whose associated $\OX$-ideal~$\widetilde \fC$ is $\mathcal C$ 
(see \cite[(2.7.11)(ii)]{EGII}). With $i\colon X_1\hookrightarrow X$ the inclusion
and $\pi\colon \OX\to i_*\sO_{\<\<X_1}$ the natural map, let $\mathcal C_1$\vspace{-.6pt} be the conductor of $\sO_{\<\<X_1}$
and let $\fC'$ be a graded $A$-ideal whose associated $\OX$-ideal  
$\<\<\widetilde{\,\<\fC'\>}\<$ is $\mathcal C'\set \pi^{-1}i_*\mathcal C_1$. Then \emph{neither\/ $\fC$ nor\/
$\fC'$ is contained in any minimal prime divisor\/ $\fP$ of\/ $(0)$ or of\/ $\fm A$.} Indeed,
such a $\fP$ is graded, and does not contain every element of positive degree in $A$
because otherwise $X_1=\proj A/\fm A$ would be empty.\vspace{.4pt}
For any homogeneous $a\notin \fP$,
$\fC_{(a)}$ is the conductor of $A_{(a)}$ (the ring of \mbox{degree-0} elements in the localization $A_a$),\vspace{-.6pt} $\fC'_{(a)}$ is the inverse image in $A_{(a)}$ of the conductor of $A_{(a)}/\fm A_{(a)}$,\vspace{.4pt}  and $\fP_{(a)}$ is a minimal prime divisor in $A_{(a)}$
of $(0)$ or of $\fm A_{(a)}$, as the case may be; and so  either $\fC\subset\fP$ or $\fC'\subset\fP$ would lead to a contradiction, as in the proof of Lemma~\ref{hexist}
(with $A$ replaced by $A_{(a)})$. Since $A_{(a)}$ and $\fm A_{(a)}$ have no embedded associated primes (see~\ref{factA}(i)), homogeneous prime avoidance implies then that there is a homogeneous $h\in \fC\cap\fC'$, of degree, say, $n>0$, such that 
$(0:_Ah)\under{1.8}{\widetilde{\phantom{iii}}}\<\<=(0)$ and $(\fm
A:_Ah)\under{1.8}{\widetilde{\phantom{iii}}}\<\<=\>\fm \OX$. Thus,
if $a\in A_{[1]}$, then the pair $(A_{(a)}, h/a^n)$ has the same
properties as the pair $(A,h)$ in Lemma~\ref{hexist}.

Globally, this $h$ gives rise to a section of $\mathcal C(n)$, i.e., to a map
$\tilde h\colon\OX\to\mathcal C(n)$. There results a composed map
$$
\fh\colon\overline{\OX\<\<}\>\>=
\overline{\OX\<\<}\>\>\otimes\OX\xrightarrow{1\otimes\tilde h\,}
\overline{\OX\<\<}\>\>\otimes\mathcal C(n)\xrightarrow{\textup{natural}\,}\OX(n).
$$
An examination of the basic definitions involved shows that for any
affine open subset $X_a\set\spec A_{(a)}\subset X$ 
($a\in A_{[1]}$), $\Gamma(X_a,\fh)$ is the composed $A_{(a)}$-homomorphism
\begin{equation}\label{haff}
\Gamma(X_a,\overline{\OX\<\<}\>\>)=
\>\overline{\<A_{(a)}\<\<}\>\>\xrightarrow{h/a^n} A_{(a)}\iso\Gamma(X_a,\OX(n))
\end{equation}
where the first map, multiplication by the regular element $h/a^n$ lying in
the conductor of~$A_{(a)}$, is injective, and the second map is the isomorphism
taking $1\in A_{(a)}$ to the section defined by $a^n\in A_{[n]}$.

We are now in a position to prove  properties of $\fh$ 
analogous to those of $h$ which were needed for the proof of Theorem~\ref{main}.

\stepcounter{slem}
\begin{slem}\label{frakhexist} 
There exist\/ $n>0$ and an\/ $\OX$-homomorphism\/ 
$\fh\colon\overline{\OX}\to\OX(n)$ such that
if\/ $R\subset T\subset K$ and $T$ is a normal  essentially-finite-type\/
$R$-algebra, or if\/ $T$ is a finite-type field extension of\/ $k,$ then 
with\/ $X_T\set X\otimes_R T,$ 
the map\/ $\fh\otimes 1$ factors uniquely as\looseness=-1 
$$
\>\overline{\<\OX\<\<}\>\>\otimes_RT
\xrightarrow[\under{1.2}{\textup{\ref{factA}}}]{\alb}
\>\overline{\<\sO_{X_{\<T}}\!}\,
\xrightarrow{\fh_{\>T}}\sO_{\<\<X_{\<T}}(n)=
\OX(n)\otimes_R T\>;
$$
and $\fh_T$ is injective, 
$\sO_{\<\<X_T}(n)/\fh_T(\sO_{\<\<X_T})$ is $T\<$-flat, and
$\sO_{\<\<X_T}(n)/\fh_T(\>\>\overline{\<\<\sO_{\<\<X_T}\<\<}\>\>)$ is
$T\<$-torsion-free.
\end{slem}

\begin{proof} The restriction of the previously described map $\fh$ to $\OX$ corresponds to multiplication by $h$ in $A$; and for any $R$-algebra $T$ at all,
the map
$$
\sO_{\<\<X_T}=\OX\otimes_R T\xrightarrow{\fh\otimes 1}\OX(n)\otimes_R T
=\sO_{\<\<X_T}(n)
$$
is given by multiplication by $h\otimes 1\in A\otimes_R T$. Let $a\in A_{[1]}$, and set $S=S^a\set A_{(a)}\otimes _R T$. As in \eqref{haff}, over $\spec A_{(a)}\subset X$ the map~$\fh$  is given by multiplication by the
regular element $h/a^n$; and similarly, over $\spec  S\subset X_T$ the map $\fh\otimes 1$ is given by multiplication 
by $h_T\set h/a^n\otimes 1$---which is regular in $S$ (see line following ~\eqref{htor}). 

Let us first check uniqueness, i.e., that\vspace{.4pt} 
if $\>\fh_T'\smcirc\alb=\fh_T''\smcirc\alb=\fh\otimes 1$, then $\fh_T'=\fh_T''$.  This can be done locally, say over  $\spec S$. Set 
$\alpha^a\set\Gamma(\spec S,\alb)$, \vspace{.6pt}  $\fh_T^{\<*\>a}\set\Gamma(\spec S,\fh_T^*)$.
\vspace{.6pt}Let $\sigma\in\overline S$, and let $g\in S$ be a regular element such that $g\sigma\in S$.
Let $\iota\colon A_{(a)}\hookrightarrow\overline{A_{(a)}}$ be the inclusion map. 
From~\ref{factA} we get that $\alpha^a\smcirc(\iota\>\otimes 1)$ is the inclusion  $S\hookrightarrow\overline S$. Then 
$$
g\>\fh_T^{\prime\>\>a}(\sigma)=\fh_T^{\prime\>\>a}(g\sigma)=
\fh_T^{\prime\>\>a}\alpha^a(\iota\>\otimes 1)(g\sigma)=
\fh_T^{\prime\prime\>a}\alpha^a(\iota\>\otimes 1)(g\sigma)=\fh_T^{\prime\prime\>a}(g\sigma)=g\>\fh_T^{\prime\prime\>a}(\sigma);
$$ 
and so, $g$ being regular,  $\fh_T^{\prime\>\>a}(\sigma)=\fh_T^{\prime\prime\>a}(\sigma)$. 
Thus $\fh_T^{\prime\>\>a}=\fh_T^{\prime\prime\>a}$, and therefore $\fh_T'=\fh_T''$.

Given the preceding uniqueness,  we need only prove existence over  \mbox{$\spec S=\spec S^a$},  since then gluing gives existence in the general case. So assume $X=\spec A_{(a)}$. When  \mbox{$R\subset T\subset K$},  \eqref{hBA}---suitably interpreted---shows that multiplication by $(h/a^n)\otimes 1$ maps
$\overline{S}$ into $S$. Call the resulting \emph{injective} map~$\fh_T^a$ (see~\eqref{haff}). The $S$-homomorphisms 
$\fh_T^a\smcirc\alpha^a\smcirc(\iota\>\otimes 1)$ and 
$\Gamma(\spec S, \fh\otimes 1)\smcirc(\iota\>\otimes 1)$ are both given by multiplication by $(h/a^n)\otimes 1$, so the restrictions of the maps $\fh_T^a\smcirc\alpha^a$ and 
$\Gamma(\spec S, \fh\otimes 1)$ to $(\iota\>\otimes 1)(S)\subset\>\>\overline{\<\<A_{(a)}\!}\,\otimes\>T$ coincide. Since  
$(h/a^n)\otimes 1$ is regular in~$S$ and 
$((h/a^n)\otimes 1)(\>\>\overline{\<\<A_{(a)}\!}\,\otimes \>T)\subset (\iota\>\otimes 1)(S)$,
therefore $\fh_T^a\smcirc\alpha^a$ and 
$\Gamma(\spec S, \fh\otimes 1)$ coincide on $\overline{\<\<A_{(a)}\!}\,\otimes \>T$,
whence $(\spec\fh_T^a)\smcirc\alb=\fh\otimes 1$, giving the asserted existence.

 A similar argument holds when
$T$ is a field extension, necessarily separable, of~$k$, except that now \eqref{hBA} 
(with $S$ in place of $A_T$) holds for $T=k$ because $h\in\fC'$, and then for any field extension $T\supset  k$ because 
$\overline{S}=(\>\>\overline{\<\<A_{(a)}\otimes_Rk\<}\>\>)\otimes_k T$ \cite[(6.14.2)]{EGA}. 

The \emph{flatness} assertion is  local on $X$, and so in view of \eqref{haff}, 
it is given by the fact that
$S/((h/a^n)\otimes 1)S=(A_{(a)}/(h/a^n)A_{(a)})\otimes_RT$ is $T$-flat (see third line before ~\eqref{htor}). Similarly, the \emph{torsion-freeness} assertion is given by Lemma~\ref{torsionfree}.
\end{proof}

Noting that with $n$ as in~\ref{frakhexist},  $\sH_{y^{}_i}\<\<(\sO_{\<\<X_i}(n))=\sH_R(\OX(n))$, we redefine $g\colon\bG\to \spec R$ to be the \emph{Quot-scheme} of $R$-flat
quotients of $\>\OX(n)$ with Hilbert polynomial
$$
\sH(n)\set\sH_{y^{}_1}\<\<\Big(\sO_{\<\<X_{\<1}}(n)\Big)-\sH_{y^{}_1}\<\<(\>\>\overline{\<\<\sO_{\<\<X_{\<1}}\<\<}\>\>)=
\sH_{y^{}_0}\<\<\Big(\sO_{\<\<X_0}(n)\Big)-\sH_{y^{}_0}\<\<(\>\>\overline{\<\<\sO_{\<\<X_0}\<\<}\>\>),
$$
see \cite[p.\,221-12, Thm.\,3.2]{G}.  For $R$-algebras $T$ and $X_T\set X\otimes_R T$ there is a
functorial bijection between $R$-morphisms $T\to\bG$ and coherent
$\sO_{\<\<X_T}$-submodules $L$ of $\sO_{\<\<X_T}(n)$ such that 
$\sO_{\<\<X_T}(n)/L$ is
$T$-flat and $\sH_T(\sO_{\<\<X_T}(n)/L)=\sH$.
The map $g$ is projective.

Assume further that $R\subset T\subset K$ and that $T$ is normal and essentially
of finite type over~$R$.  
 Let $\psi\colon\spec K\to\bG$ be the $R$-morphism corresponding to 
$\fh_K(\>\>\overline{\<\<\sO_{\<\<X_0}\<\<}\>\>)
\subset\sO_{\<\<X_0}(n)$.  
Suppose  the domain of definition $\textup{D}(\psi_T^{})$\vspace{.4pt}
of the corresponding rational $R$-morphism $\psi_T^{}\colon T\to \bG$
is all of~$\spec T$.  Then $\psi_T^{}$
determines a coherent $\sO_{\<\<X_T}$-module
$L_1\subset\sO_{\<\<X_T}(n)$ with $\sO_{\<\<X_T}(n)/L_1$
\vspace{.6pt} $T$-flat, such that $L_1\otimes_T K
= \fh_K(\>\>\overline{\<\<\sO_{\<\<X_0}\<\<}\>\>)=L_2\otimes_T K$ where $L_2\set\fh_T(\>\>\overline{\<\<\sO_{\<\<X_T}\<\<}\>\>)$;\vspace{.6pt} and 
by~\ref{frakhexist},  $\sO_{\<\<X_T}(n)/L_2$ is
$T$-torsion-free. The assertion being local, Lemma~\ref{uniquetorsion}
shows then that $L_1=L_2$, so that\vspace{.4pt}
$\sO_{\<\<X_T}(n)/L_2\cong \sO_{\<\<X_T}(n)/L_1$ \emph{is $T$-flat,} and
there is an exact sequence of $T$-flat $\sO_{\<\<X_T}$-modules
$$
0\to \>\>\overline{\<\<\sO_{\<\<X_T}\<\<}\>\>\xrightarrow{\fh_T}\sO_{\<\<X_T}(n)\to 
\sO_{\<\<X_T}(n)/L_2\to 0.
$$

If $T=R$, so that $\sO_{\<\<X_T}=\OX$,  then we deduce that 
$\fh\otimes 1\colon\>\>\overline{\<\<\OX\<\<}\>\>\otimes_R k\to
\sO_{\<\<X_{\<1}}\<\<(n)$ is injective, whence, by~\ref{frakhexist}, so is 
$\alb\colon\>\>\overline{\<\<\OX\<\<}\>\>\otimes_R k\to
\>\overline{\<\sO_{X_{\<1}}\!}\,\>$---in other words, the 
schematic dominance mentioned in the first paragraph of the proof of Theorem~\ref{main2} holds.
\vspace{1.5pt}

Thus it suffices to show that \emph{if\/
$\sH_{y^{}_1}\<\<(\>\>\overline{\<\<\sO_{\<\<X_1}\<\<}\>\>)=
\sH_{y^{}_0}\<\<(\>\>\overline{\<\<\sO_{\<\<X_0}\<\<}\>\>)$ 
then\/ $\textup{D}(\psi)$ is all of}
$\>\spec R$.\vspace{1.5pt}

To do so we proceed in essentially the same way as in the corresponding
part of the proof of Theorem~\ref{main}. Let $(S,\fm_S,k_S)$ be as in that proof.  The natural map $\spec S\to \bG$, which localizes generically to $\psi$ and hence is the same as the above $\psi_S$, corresponds to 
$\fh_S(\>\>\overline{\<\<\sO_{\<\<X_S}\<\<}\>\>)$, which is $S$-flat 
(so that $\overline{\<\<\sO_{\<\<X_S}\<\<}\>\>$  is $S$-flat); and
\begin{align*}
\sH_{\fm_S}(\>\>\overline{\<\<\sO_{\<\<X_S}\<\<}\>\>\otimes_S k_S)=
\sH_S(\>\>\overline{\<\<\sO_{\<\<X_S}\<\<}\>\>)
&=
\sH_S\Big(\sO_S(n)\Big)-\sH_S(\sO_S(n)/\fh_S\Big(\>\>\overline{\<\<\sO_{\<\<X_S}\<\<}\>\>)\Big)\\
&=\sH_{y_0}\Big(\sO_{\<\<X_0}(n)\Big)-\sH(n)
=\sH_{y^{}_0}\<\<(\>\>\overline{\<\<\sO_{\<\<X_{0}}\<\<}\>\>)
=\sH_{y^{}_1}\<\<(\>\>\overline{\<\<\sO_{\<\<X_{\<1}}\<\<}\>\>).
\end{align*}
 Since $k_S$ is separable over $k$, therefore (by \cite[(6.14.2)]{EGA}),
 $$
 \sH_{\fm_S}(\>\>\overline{\<\<\sO_{\<\<X_{k_S}}\<\<}\>\>)=
  \sH_{\fm_S}(\>\>\overline{\<\<\sO_{\<\<X_1}\<\<\otimes_k k_S}\>\>)=
 \sH_{\fm_S}(\>\>\overline{\<\<\sO_{\<\<X_1}\<\<}\>\>\otimes_k k_S) =
  \sH_{y^{}_1}\<\<(\>\>\overline{\<\<\sO_{\<\<X_1}\<\<}\>\>).
 $$
 where the last equality holds\vspace{.4pt} because cohomology is compatible with flat base change \cite[(1.4.15)]{EGIII}.
Moreover, by ~\ref{injectivity}, the  map\vspace{.8pt} 
$\>\>\overline{\<\<\sO_{\<\<X_{\<S}}\<\<}\>\>\otimes_S k_S\to
\>\>\overline{\<\<\sO_{\<\<X_{k_{\<S}}}\<\<}\>\>$ in~\ref{frakhexist} (applied to $f\otimes 1\colon X_S\to \spec S$, with $T=k_S$)
is injective.\vspace{.4pt} Since the source and target of this map have the same Hilbert polynomial,
it must be an isomorphism, and so
$$
\fh_{k_{\<S}}(\>\>\overline{\<\<\sO_{\<\<X_{k_{\<S}}}\<\<}\>\>)=
(\fh_S\otimes 1)(\>\>\overline{\<\<\sO_{\<\<X_S}\<\<}\>\>\otimes_S k_S)
\subset\sO_{\<\<k_S}(n).
$$
This signifies that the natural composition\vspace{.6pt} $\spec k_S\hookrightarrow\spec S\to \bG$ corresponds to
\mbox{$\fh_{k_S}(\>\>\overline{\<\<\sO_{\<\<X_{k_S}}\<\<}\>\>)\subset\sO_{\<\<k_S}(n)$.}
But $\fh_{k_S}(\>\>\overline{\<\<\sO_{\<\<X_{k_S}}\<\<}\>\>)
=(\fh_k\otimes 1)(\>\>\overline{\<\<\sO_{\<\<X_1}\<\<}\>\>\otimes_k k_S)$.

Thus if $\psi_1\colon \spec k\to\bG$ corresponds to 
$\fh_k\>\>\overline{\<\<\sO_{\<\<X_1}\<\<}\>\>\subset \sO_{\<\<X_1}(n)$ then the following 
natural diagram commutes:
$$
\begin{CD}
\spec k_S @>>>\spec S \\
@VVV @VVV \\
\spec k @>>{\psi_1}>\bG
\end{CD}
$$
As in the proof of Theorem~\ref{main}, the desired conclusion results. 
\end{proof}

\section{Application to complex spaces.}

In this section we translate Theorems~\ref{main} and~\ref{main2} from
algebraic into analytic geometry, that is, from the context of schemes
to the context of complex spaces, see Theorems~\ref{mainan}
and~\ref{main2an}.  The treatment of these theorems comes after a few
pages of preliminaries on equinormalization for maps of complex
spaces.

\pagebreak[3]

To establish terminology and notation, referring to
\cite[expos\'es 9,\ 10,\ 13]{C} or \cite{Gr} for details we first review
a few foundational facts about complex spaces, beginning with their
definition. For any $n\ge0$ and open $U\subset \C^n$ let $\sO_U$ be
the sheaf of holomorphic functions $U\to \C$ (functions given locally
by convergent power series).  The \emph{local model} determined by an
open $U\subset \C^n\ (n\ge0)$ and a finite-type $\sO_U$-ideal~$\cI$,
is the topological space $X_\cI\set{}$support of $\sO_U\</\cI$
together with the sheaf of $\C$-algebras $(\sO_U\</\cI)|_{X_\cI}$. A
\emph{complex space} is a $\C$-ringed space---a topological space~$X$
with a sheaf~$\>\OX$ of $\C$-algebras---such that each point of~$X$
has an open neighborhood isomorphic as a $\C$-ringed space to some
local model. \looseness=-1

Any open subset $U$ of a complex space~$X$ can, and will, be viewed as
a complex space, with $\sO_U\set\OX|_U$. A subset of $X$ will be
called \emph{analytic} if it is the underlying topological space of a
closed complex subspace of~$X$, i.e., if it is the support of
$\OX/\cI$ for some finite-type $\OX$-ideal~$\cI$.

An \emph{analytic} (or \emph{holomorphic}) \emph{map} 
$f\colon(X,\OX)\to (Y,\sO_Y)$ \emph{of complex spaces} is a
continuous map $f_0\colon X\to Y$ together with a homomorphism of sheaves of
$\C$-algebras $f_1\colon \sO_Y\to f_{0*}\OX$.  The sheaf $\>\OX$ is
canonically isomorphic to the sheaf of (germs of) analytic\- maps from
$(X,\OX)$ to~$(\C^1\<,\>\sO_{\C^1}\<)$, and $f_1$ is given locally by
composition of analytic maps. It is a consequence of (indeed,
essentially equivalent to) R\"uckert's nullstellensatz that if $X$ is
\emph{reduced}, i.e., $\OX$ has no nonzero nilpotents, then $f$ is
uniquely determined by~$f_0$ \cite[p.19-18, Cor.5]{C}; in particular,
$\OX$ is canonically isomorphic to a subsheaf of the sheaf of
continuous $\C$-valued functions on~$X$.

Though principally interested here in reduced complex spaces, we must
allow for~ the possible occurrence of nonreduced ones, for instance as
fibers of morphisms $f\colon X\to Y$ of reduced complex spaces. (The
fiber $X_y\set f^{-1}y\ (y\in Y)$ is by definition the complex space
$X\times_Y[\>y]$ where $[\>y]$ is the complex subspace
$(\{y\},\C)\subset(Y,\sO_Y)$.)\vspace{1pt}

A \emph{pointed complex space} is a pair $(X,x)$ where $X$ is a
complex space and $x\in X$. A~\emph{morphism $(X,x)\to (Y,y)$ of
pointed complex spaces,} or \emph{pointed analytic map,} is an
analytic map $f\colon X\to Y$ such that $f(x)=y$. The localization of
the category of pointed complex spaces with respect to open immersions
is the category $\mathsf{AG}$ of complex germs. Explicitly, the
objects of $\mathsf{AG}$---called \emph{germs}---are the pointed
complex spaces, and the morphisms $(X,x)\to (Y,y)$ in
$\mathsf{AG}$---called \emph{map germs}---are equivalence classes of
pointed analytic maps $(U,x)\to (Y,y)$ with domain $U$ an open
neighborhood of $x$ in $X\<$, where two such maps are equivalent if they
agree on an open neighborhood of $x$ contained in both their
domains. The composition of two map germs is defined in the obvious
way. Two germs $(X,x)$ and $(Y,y)$ are isomorphic iff there exist open
neighborhoods $U$ and~$V$ of $x$ and~$y$ respectively, such that the
pointed spaces $(U,x)$ and $(V,y)$ are isomorphic.

\enlargethispage*{7pt}
The category $\mathsf{AA}$ of \emph{analytic $\C$-algebras} has as its
objects homomorphic images of rings of convergent power-series with
coefficients in $\C$, and as its morphisms all $\C$-algebra
homomorphisms. Analytic $\C$-algebras are \emph{excellent} noetherian
local rings (see footnote under Remark~\ref{spadelem}(a)), and
$\C$-algebra homomorphisms between them are automatically local (i.e.,
take nonunits to nonunits). 

Associating to any germ the stalk
$\sO_{\<\<X\<,\>x}$, and to any map germ $(X,x)\to (Y,y)$ the induced
map $\sO_{Y,\>y}\to \sO_{\<\<X\<,\>x}$ one gets a contravariant
functor $\sO\colon\mathsf{AG}\to\mathsf{AA}$ which is in fact an
\emph{antiequivalence of categories} \cite[p.\,13-2,
Thm.\,1.3]{C}. This
is a fundamental link between local algebra and local
analytic geometry. The connection between algebra and geometry is
illustrated in the following Proposition.

\begin{defn}\label{condan}
 A map of complex spaces satisfies $\condan$ if
\begin{itemize} \item $f$ is flat\/ \vspace{-.4pt}$($i.e., the induced
map $\sO_{Y,\>f(x)}\to \sO_{\<\<X\<,\>x}$ is flat for all\/ $x\in X)$
and all the fibers of $f$ are reduced.
\item $Y$ is normal\/  $($i.e.,  $\sO_{Y,\>y}$ is normal for all\/ $y\in Y)$.
\item $X$ is locally equidimensional\/ $($i.e., for each\/ $x\in X$
the irreducible components of\/~$X$ passing through\/ $x$ all have the
same dimension, i.e., the connected components of $X$ are
equidimensional). \emph{Equivalently} (see \cite[p.\,116, 15.1]{Ma}), 
the fibers of $f$ are locally equidimensional.
\end{itemize}
\end{defn}

With $y\set f(x)$ and $\fm_y$ the maximal ideal of~$\sO_{Y\!,\>\>y}$,
the local ring of $x$ on the analytic fiber $f^{-1}y$ is 
$\sO_{\<\<X\<,\>x}/\fm_y\sO_{\<\<X\<,\>x}$ \cite[p.\,10-07]{C}. It
follows from ~\ref{Nishimura} and \cite[p.\,184, Cor.\,(ii)]{Ma}
that \emph{if\/ $f$ satisfies~$\condan$ then\/ $X$ is reduced.}

\begin{prop} \label{cluban}
The map\/ $\spec\phi\colon\spec R\to\spec A$ associated to a\/
$\C$-homomorphism of analytic\/ $\C$-algebras\/ $\phi\colon R\to A$
satisfies\/ $\cond$ iff\/ $\phi$ is isomorphic to\/ $\sO(F)$ where\/
$F$ is a map germ represented by a pointed analytic map\/ $f\colon
(X,x)\to (Y,y)$ satisfying\/ $\condan$.

Moreover, if\/ $\spec\phi$ does satisfy\/ $\cond$ then\/ such an\/ $f$ exists 
with\/ $X$  reduced and connected,
$Y$ normal and irreducible, $f$ surjective, and all
fibers of\/ $f$ of pure dimension\/ $\dim_x X-\dim_y Y;$ and if, in addition, 
$A/\fm A$ is normal, then there is an\/ $f$ with the preceding properties and
also such that\/ $X$ is normal and all the fibers of\/ $f$ are normal. 

%Moreover, with\/ $\fm$ the maximal ideal of\/ $R,$ $A/\fm A$ is normal
%iff\/ $f$ can be chosen so that all its fibers are normal---and in
%this case one can assume that\/ $X$ is normal and irreducible.

\end{prop}

\begin{proof}
Suppose given $f$ satisfying $\condan$ 
%(resp.~and in addition having normal fibers), 
and $\phi\cong\sO(F)$. The closed fiber of
 $\spec\phi$ is isomorphic to the $\spec$ of the local ring of $x$ on
 the analytic fiber $f^{-1}y$, see \cite[p.\,10-07]{C}; hence
 $\spec\phi$ is flat, with geometrically reduced 
 %(resp.~normal) 
 closed fiber, see~\ref{geornr}.  The local ring $R\cong\sO_{Y,\>y}$ is
 normal, and, being excellent, satisfies $\ff{nor}$. So by
 ~\ref{Nishimura}, all the fibers of $\spec\phi$ are geometrically
 reduced, 
% (resp.~normal), 
 so that $\spec\phi$ is reduced.
 %(resp.~normal). 
 The local ring $A\cong\sO_{\<\<X,\>x}$, being
 excellent, is a universally catenary Nagata ring, whose
 equidimensionality results from \cite[p.\,20-03, Cor.\,2]{C}. Thus
 $\spec\phi$ satisfies~\cond.
% ; and since $A/\fm A$ is the local ring of a point on a fiber 
% of~$f$, if the fibers of~$f$ are normal then $A/\fm A$ is (geometrically) normal.

Conversely, suppose that $\spec\phi$ satisfies~\cond\ (resp.~and in addition $A/\fm A$ is normal).
By the above-mentioned equivalence of categories, $\spec\phi\cong \sO(F)$ for some map germ~$F$ represented by a pointed analytic map 
$f\colon (X,x)\to (Y,y)$. For any open neighborhoods $U$ of~$x$ and $V$ of~$y$ with  $f(U)\subset V$, we are free to replace $f$ by \mbox{$f|_U\colon(U,x)\to (V,y)$,}
a process referred to as ``shrinking." 

Since $A$ is equidimensional, \cite[Chap.\,5, \S4.2, p.\,106]{Gr}
shows that after shrinking $X$ if necessary, we may assume that all
components of $X$ have dimension~$d\set\dim_x X$ at each of their
points. Further, Frisch's theorem \cite{Fr} allows to assume, after
more shrinking, that $f$ is flat. Still further, \cite[p.\,160,
3.22]{Fi} allows us to assume that all the fibers of~$f$ are reduced
(resp.~normal).  Then by ~\ref{Nishimura} and \cite[p.\,184,
23.9]{Ma}, $X$ is reduced (resp.~normal).  Since $y$ is a normal point
(because $\sO_{Y,\>y}\cong R$) and the nonnormal points of $Y$ form an
analytic subset \cite[p.\,21-09, Thm.\,2]{C}, $Y$ has a normal open
neighborhood $V\<$, which may be assumed to be connected, hence
irreducible and of dimension $e\set\dim_yY$ at each of its points.  It
follows from the dimension relation for flat local homomorphisms
\cite[p.\,116, 15.1]{Ma}, applied to
$\sO_{Y\<,\>\>f(s)}\to\sO_{\<\<X,\>s}\ (s\in f^{-1}V)$, that after $Y$
is shrunk to $V$ and $X$ to $f^{-1}V$, every fiber of $f$ has pure
dimension $d-e$. Since flat maps are open, \cite[p.\,156, 3.19]{Fi},
therefore $f$ maps the connected component $X_0$ of $X$ containing $x$
onto a connected open subset $Y_0\subset Y$. Shrinking $X$ to $X_0$
and $Y$ to $Y_0$ finishes the proof.\vspace{-2pt}
\end{proof}

Equinormalizability of a \emph{reduced} map 
$f\colon X\to Y$ of analytic spaces (i.e., $f$~is flat, with reduced
fibers) means, again, the existence of a simultaneous
normalization, now defined as follows. Recall first that an analytic
map $\mu\colon Z'\to Z$ is \emph{finite} if it is proper and has
finite fibers; that $\mu$ is \emph{bimeromorphic} if there are
analytically rare (hence {nowhere dense}, and  conversely when $X$ is reduced \cite[pp.\,38--40, \S0.43]{Fi})
%
%\footnote{\emph{Analytically rare}  implies \emph{nowhere dense}, and  conversely when $X$ is reduced \cite[pp.\,38--40, \S0.43]{Fi}.}
analytic subsets $W\subset Z$ and $W'= \mu^{-1}(W)\subset Z'$ such
that $\mu$ induces an isomorphism from $Z'\setminus W'$ onto
$Z\setminus W$; and that $\mu$ is a \emph{normalization map,} or a
\emph{normalization of\/~$Z,$} if $\mu$ is finite and bimeromorphic,
and in addition $Z$~is reduced and $Z'$~is normal. For fixed $Z$, any
two normalizations are isomorphic \cite[p.\,21-11, Cor.\,3]{C}.  We
will use freely properties of normalization found, e.g., in
\cite[Chapter 8]{Gr} or \cite[Expos\'e 21, \S4]{C}.\looseness=-1

\begin{defn} A \emph{simultaneous normalization} of a 
reduced analytic map $f\colon X\to Y$  is a finite analytic map 
$\nu\colon Z\to X$ such that such that 
$\bar f\set f\<\smcirc\nu$ is \emph{normal} (flat, with all nonempty
fibers normal), and such that for each $y\in f(X)$ the induced map of
fibers $\nu_y\colon \bar f^{-1}y\to f^{-1}y$ is a normalization map.
\end{defn}

\begin{prop} \label{simnorm}
Let\/ $f\colon X\to Y$ be an analytic map satisfying\/ $\condan$. Then$\>:$
\vspace{1pt}

{\rm(i)} Any simultaneous normalization of\/ $f$ is a normalization of\/ $X\<$.

{\rm(ii)} A normalization\/ \mbox{$\nu\colon \BX\to X$} is a simultaneous normalization of\/ $f$ if and only if the nonempty fibers of $f\nu$ are all normal.
\end{prop}

\begin{proof}
It suffices to prove the assertions for the restriction of $f$ to any of the connected components of~$X$. So we may assume that $X$ and $Y$ are connected. Then $Y$, being normal, is irreducible, and so has the same dimension, say $e$, at each of its points. Moreover, since $X$ is locally equidimensional, for any 
$n\ge 0$ the (locally finite) union of the $n$-dimensional irreducible components of~$X$ is an open and closed complex subspace; and hence we may assume that $X$ has the same dimension, say $d$, at each of its points. 

(i) It follows from~\ref{cluban} that every  fiber 
$f^{-1}y\ (y\in f(X))$ has pure dimension $d-e\>$. For a simultaneous normalization 
$\nu\colon Z\to X$ the same is therefore true of~$(f\nu)^{-1}y$, 
the normalization of~$f^{-1}y$
(because a proper bimeromorphic map  $\mu\colon W'\to W$ maps each irreducible component  $C\subset W'$ bimeromorphically onto an irreducible component of $W$ having the same dimension as~$C$.)
The restriction of $f\nu$ to any open subset $U$ of $Z$ meeting only one irreducible component $C$ of $Z$ is flat, so \cite[p.\,116, 15.1]{Ma} implies that $\dim C=\dim U=e+(d-e)=d$. Since $\nu$ is finite, we conclude that
\emph{$\nu(C)$~is an irreducible component of\/ $X\<$.}

By~\ref{Nishimura} and \cite[p.\,184, Cor.\,(ii)]{Ma}, $X$ is reduced and $Z$ is normal. It remains then to show that  $\nu$ is bimeromorphic. 
\begin{slem}
For any\/ $x\in X$ the following conditions are equivalent.

{\rm(a)} $x$ is a normal point of\/ $f^{-1}f(x)$.

{\rm(b)} There is a unique\/ $z\in Z$ such that $\nu(z)=x;$ and for that $z,$ the induced map germ $(Z,z)\to (X,x)$ is an isomorphism.
\end{slem}

\pagebreak[3]
Once the lemma is proved, then since the set of points $W\subset X$ not satisfying (a) is analytic \cite[p.\,160, Proposition]{Fi}, clearly nowhere dense 
(because $f$~has reduced fibers), we see, by (b), that the restriction of $\nu$ maps $Z\setminus \nu^{-1}(W)$ one-to-one and locally isomorphically,
hence globally isomorphically, onto  $X\setminus W\>$; and moreover, $\nu^{-1}(W)$ is nowhere dense in~$Z$, because $\nu$ maps components of~$Z$ onto components of~$X$. Thus $\nu$ is  bimeromorphic.%
\footnote{One way to avoid the referenced theorem is via differentials and ramification \cite[p.20-12, Corollaire]{C}. }

As for the proof of the lemma, the implication (b)${}\Rightarrow{}$(a) is simple, because
the map $f\nu$ has normal fibers. Conversely, if (a) holds, then, by standard properties of normalization, since $\nu$ normalizes fiberwise, there just one $z$ with $\nu(z)=x$, and with $\fm$ the maximal ideal of
$\sO_{Y,\>f(x)}$, $\nu$ induces an isomorphism 
$\sO_{\<\<X,\>x}/\fm\sO_{\<\<X,\>x}\iso\sO_{Z,\>z}/\fm\sO_{Z,\>z}$. Then,\vspace{.4pt} since $\sO_{Z,\>z}$ is a finite $\sO_{\<\<X,\>x}$-module, Nakayama's lemma shows that $\nu$ induces a surjection $\sO_{\<\<X,\>x}\twoheadrightarrow \sO_{Z,\>z}$, which must be an isomorphism, because $\sO_{\<\<X,\>x}$ is normal (see~\ref{cluban}) and, as above,  
$\dim \sO_{Z,\>z}=\dim\sO_{\<\<X,\>x}=d$. 

This completes the proof of (i).
\vspace{1pt}
 
(ii) By definition, if $\nu$ is a simultaneous normalization then the nonempty  fibers of~$f\nu$ are normal.

Suppose, conversely, that the fibers of $f\nu$ are normal.
Let $y\in Y\<$ and let $C$ be an irreducible component of~$\BX_{\!y}\set(f\nu)^{-1}y$.
Recall that $\nu(C)$ can be viewed as a
closed complex subspace of $X_y\set f^{-1}y$ \cite[p.\,65, Proposition]{Gr}. We claim that \emph{$\nu(C)$ is an irreducible component of} $X_y\set f^{-1}y$. 
To see this, let $z\in C\>$ lie on no other component of~$\BX_{\!y}$ and set $x\set\nu(z)$. Since $\nu$ is a normalization and $X$ is equidimensional, therefore $\dim_z\<\BX=\dim_x X\<$. Also, 
by \cite[p.\,116, 15.1(ii)]{Ma},
$$
\dim_x\nu(C)=\dim_z C=
\dim_z\<\BX_{\!y}\ge\dim_z\< \<\BX-\dim_y Y=\dim_x X-\dim_y Y = \dim_x X_y,
$$
and the assertion results.

Next, let $D$ be the dense open subset of normal points of the reduced analytic 
space~$X_y$. By~\ref{Nishimura} and \cite[p.\,184, Cor.\,(ii)]{Ma}, 
every point $w\in D$ is normal on $X\<$, and so $w$ has an open $X$-neighborhood~$U_w$ such that 
$\nu^{-1}(U_w)\to U_w$ is an isomorphism. Every component~$C$ of~$\BX_{\!y}$ must meet $\nu^{-1}(U_w)$ for some $w$, 
since $D$ meets the component $\nu(C)$ of~$X_y$ in an open dense set. 
It follows then from  \cite[p.\,21-11, Cor.\,3]{C} that  $\nu\colon \BX_{\!y}\to X_y$ is a normalization; and by~\ref{birational}(iv), $f\nu$ is flat. Thus  $\nu$ is a simultaneous normalization of~$f$.\looseness=-1
\end{proof}

We say that an analytic map $f\colon X\to Y$ is \emph{equinormalizable at $x\in X$} if
the induced scheme-map $\spec\sO_{\<\<X\<,\>x}\to\spec\sO_{Y\<\<,\>\>f(x)}$ is equinormalizable.

\begin{scor}\label{equinorm local} Let\/ $f\colon X\to Y$ be an analytic map satisfying\/ $\condan$.\vspace{1pt}

{\rm(i)} If\/ $f$ is equinormalizable at\/ $x\in X$ then\/ the restriction of $f$ to some neighborhood of\/~$x$ is equinormalizable.\vspace{1pt}

{\rm(ii)} $f$ is equinormalizable iff\/ $f$ is equinormalizable at\/ $x$ for all\/ $x\in X$.
\end{scor}

\begin{proof}
 Let $\nu\colon Z\to X$ be a normalization, and let $z_1,\dots, z_n$ be the distinct
points in~$\nu^{-1}x$.  The integral closure of $\sO_{\<\<X,\>x}$ is 
$B\set\prod_{i=1}^n\sO_{Z, z_i}$, which is equidimensional since $\sO_{\<\<X,\>x}$ is.

Let $(R,\fm)$ be the local ring $\sO_{Y,\>f(x)}$. 
By~\ref{Snormal}, a simultaneous normalization of schemes is a normalization. So equinormalizability of $f$ at $x$ implies that $B$ is $R$-flat and that $B/\fm B$ is normal, whence, by ~\ref{Nishimura}, the natural map $\spec B\to\spec R$ satisfies $\cond$. By~\ref{cluban}, $f\nu$ restricted to some neighborhood~$V$ of $\nu^{-1}x$ has normal fibers. Since $\nu$ is finite 
there is a neighborhood $W$ of $x$ such that $\nu^{-1}W\subset V$.  By~\ref{simnorm}(ii), the restriction of $f$ to $W$ is equinormalizable, proving~(i).\looseness=-1

As for (ii), if $f$ is equinormalizable at every $x\in X$ then it results from (i) that the nonempty fibers of $f\nu$ are normal, so that $f$ is equinormalizable.

Conversely, if $f$ is equinormalizable, so that the nonempty fibers of $f\nu$ are normal, then by~\ref{cluban}, \ref{Nishimura} and~\ref{CorSnormal},
$f$ is equinormalizable at every $x\in X$.
\end{proof}

\begin{scor}\label{fibers}
If the set\/ $S_{\<\<f}\set\{\,x\in X\mid x\textup{ not normal on }X_{\<\<f(x)}\,\}$ 
$($which is analytic\/ \textup{\cite[p.\,160, 3.22]{Fi}}$)$ is proper over\/ $Y$ then
there is a nowhere dense analytic subset\/ $T\subset Y$ such that\/ $f$ is equinormalizable at every\/ $x\notin f^{-1}\>T$.
\end{scor}

\begin{proof}
Let $\nu\colon\BX\to X$ be a normalization.  Set
$S_{\<f\nu}\set\{\,z\in\BX\mid z\textup{ not normal on }
\BX\!_{f\nu(z)}\,\}$, an analytic subset of~$\BX$.  Then 
$\nu(S_{\<f\nu})\subset S_{\<\<f}$---for if $x=\nu(z)$ is normal on 
$X_{\<\<f(x)}$ then $x$ is normal on $X$ (see~\ref{cluban}), 
so that for some neighborhood $U$ of $x$, $\nu$ induces an isomorphism 
$\nu^{-1}U\iso U\<$, whence $z\in\nu^{-1}U$ is normal on  
$\BX\!_{f\nu(z)}$. By ~\ref{simnorm}(ii), $f$ is equinormalizable at
every point outside $\nu(S_{\<f\nu})$. Since $S_{\<\<f}$ is proper
over~$Y\<$, the set $T\set f\nu(S_{\<f\nu})$, which is analytic
\cite[p.\,213]{Gr}, and nowhere dense in $Y$ \cite[p.\,270,
Thm.\,(I.5)]{Mn}, is as asserted.
\end{proof}
%:curves

From now on, by \emph{curve} we mean \emph{reduced\/ $one$-dimensional
complex space with no isolated points.}  An \emph{analytic family of
curves} is a flat map $f\colon X\to Y$ of complex spaces whose fibers
$X_y\ (y\in Y)$ are curves. If $Y$ is normal 
then such an $f$ satisfies $\condan$, and $X$ is reduced.
 For families of curves with $Y$ reduced 
a simultaneous normalization is the same thing as 
a ``very weak simultaneous resolution" \cite[p.\,72, D\'efn.\,2]{T2}. 

The following theorem~\ref{mainan} gives a numerical criterion of
equinormalizability at a point $x\in X$ for a family of curves
$f\colon X\to Y\<$.  Arguing as in~\ref{cluban}, and with the aid of
the dimension relation \cite[p.\,116, 15.1]{Ma}, we see that a pointed
analytic map $(X,x)\to(Y,y)$ defines a map germ represented by a
family of curves $f\colon X\to Y$ with normal $Y$ if and only if
$(A,\fM)\set\sO_{\<\<X,\>x}$ is flat over $(R,\fm)\set\sO_{Y\<,\>y}$,
$A/\fm A$ is reduced and of pure dimension~1, and $R$ is normal. For
such a family the set Sing$(f)$ of points on $X$ where the fiber
$f^{-1}f(x)$ is singular is analytic \cite[p.\,160, 3.22]{Fi}, meeting
each fiber in a discrete set; so we may, and will, assume that the
induced map Sing$(f)\to Y$ is \emph{finite}, and that Sing$(f)\cap
f^{-1}y=\{x\}$.\vspace{1pt}

The theorem is due to Teissier when $Y=\C^1$, and was extended to arbitrary normal~$Y$ by Raynaud (see remarks in the Introduction).  Here we will deduce the theorem from its algebraic counterpart Theorem~\ref{main}.

 \begin{defn}\label{geodelta}
For a curve $C$ and $x\in C$,  set
$$
\dt(C,x)\set
\dim_{\>\C}(\>\overline{\<\sO_{\<\<C,\>x}\<\<}\>\>/\sO_{\<\<C,\>x}).
$$
If $C$ has only finitely many singular points, we set
$$
   \dt(C)\set \underset{x\in C}{\sum}\,\dt(C,x).
$$
\end{defn}

\begin{thm}\label{mainan}
 Let\/ $f\colon (X,x_0^{})\to (Y,y^{}_0)$ be a flat pointed analytic map, with\/ $Y$ normal, such that the nonempty fibers of\/ $f$ are curves.  Assume that\/ $\textup{Sing}(f)$ is finite over\/~$Y$ and that $\textup{Sing}(f)\cap f^{-1}y^{}_0 =\{x_0^{}\}$. Then\/ $f$ is equinormalizable at\/ $x_0^{}$ if and only if\/ $\dt(X_y)$ is  the same for all~$y$ in some neighborhood of\/~$y^{}_0$.
 \end{thm}

\begin{proof} As above, $f$ satisfies $\condan$ and
 $X$ is reduced.
Since flat analytic maps are open
\cite[p.\,156, 3.19]{Fi}, we may assume further that $f$ is surjective.

Let $\nu\colon\BX\to X$ be a normalization, and set $\overline{\OX}\set\nu_*\sO_{\<\<\BX}$.
As in~\ref{cluban} we find that if $x\in X$ is a normal point of $f^{-1}f(x)$ then $x$ is normal on~$X$. Hence the support of the coherent sheaf $\overline{\OX}/\OX$---topologically the set of nonnormal points of $X$---is finite over~$Y\<$, and the 
$\sO_Y$-module $f_*(\>\>\overline{\<\<\OX\<}/\OX)$ is \emph{coherent} \cite[p.\,64, {\bf3}]{Gr}. After shrinking we may assume that $Y$~is irreducible and that $f_*(\>\>\overline{\<\<\OX\<}/\OX)$ is the cokernel of a map $\sO_Y^m\to\sO_Y^n$ represented by an $n\times m$ matrix whose entries are holomorphic functions on $Y$. The rank~$r_y$
of this matrix when evaluated at the point $y\in Y$ is a function on $Y$ taking its maximum value~$r$ on a dense open subset $U\subset Y\<$, the complement of the analytic set defined by the vanishing of the $r\times r$ subdeterminants. For any 
$y\in Y$ the germs at $\sO_y\set\sO_{Y\<,\>y}$ of these $r\times r$ subdeterminants are not all zero, since every neighborhood of~$y$ meets $U$; while the germs of the $(r+1)\times(r+1)$ subdeterminants do vanish. Setting $K_y\set{}$the fraction field of $\sO_y$, we deduce that \emph{the dimension of the\/ $K_y$-vector space\/ 
$f_*(\>\>\overline{\<\<\OX\<}/\OX)_y\otimes_{\sO_y} K_y$ is\/ $(n-r)$}---which does not depend on\/  $y$.

Since for any $y\in Y$, $f_*(\>\>\overline{\<\<\OX\<}/\OX)_y=\oplus_{x\in f^{-1}y}(\>\>\overline{\<\<\OX\<}/\OX)_x$ (where $x$ contributes nothing to the sum unless $x\in N$, the nonnormal locus of~$X$), we have 
\begin{equation}\label{gendtsum}
n-r=\sum_{x\in N\cap f^{-1}y} \dim_{K_y}\! \Big((\overline{\OX}/\OX)_x \otimes_{\sO_y} K_y\Big).
\end{equation}

\stepcounter{slem}

\begin{slem}\label{dtsum}
If the fiber\/ $\BX_{\!y}\set(f\nu)^{-1}y$ is normal then\/ $x\in X_y$ is  normal 
iff\/ $x$ is normal  on\/ $X$. Thus, 
$$
\dt(X_y)=\sum_{x\in N\cap f^{-1}y}\dt(X_y,x)=
\sum_{x\in N\cap f^{-1}y}\dim_{\>\C}\! \Big(\>\>\overline{(\<\<\OX_{\<\<y})_x}/(\OX_{\<\<y})_x\Big).
$$
\end{slem}

\begin{proof}
That [$x$ normal on $X_y$] implies [$x$ normal on~$X$] is shown as in ~\ref{cluban}.
Conversely, if $x$ is normal on $X$, then $x$ has a normal neighborhood $V\<$, and the normalization~$\nu$ maps~$\nu^{-1}V$ isomorphically onto $V\<$; 
so $X_y\cap V\cong \BX_{\!y}\cap\nu^{-1}V$ is normal.
 \end{proof}

Now suppose $f$  equinormalizable at $x_0^{}$. By~\ref{equinorm local} we may assume (after shrinking~$X$) that $f$ is equinormalizable for all $x\in X$. Then  ~\ref{birational}(iii)$'$ holds for 
$\mu=\nu$ (see~\ref{CorSnormal}),  and hence, with $k_y\cong\C$ the residue field of $\sO_{Y\!,\>\>y}$, 
$(\>\>\overline{\<\<\OX\<}/\OX)_x$ is a free $\sO_y$-module, of rank 
$$
\dim_{K_y}\! \Big((\>\>\overline{\<\<\OX\<}/\OX)_x \otimes_{\sO_y} K_y\Big)=
\dim_{\>\C}\! \Big((\>\>\overline{\<\<\OX\<}/\OX)_x \otimes_{\sO_y} k_y\Big)=
\dim_{\>\C}\! \Big(\>\>\overline{(\<\<\OX_{\<\<y})_x}/(\OX_{\<\<y})_x\Big),
$$
where the last equality holds because by the definition of simultaneous normalization, 
$\>\>\overline{(\<\<\OX_{\<\<y})_x}=\overline{\sO_{\<\<X,\>x}}\otimes_{\sO_y} k_y$. So by~\eqref{gendtsum} and~\ref{dtsum}, $\dt(X_y)=n-r$ for all $y$.

\vspace{1.5pt} 
Assume conversely that $\dt(X_y)$ is constant near $y^{}_0$.  There exists
a sequence $y^{}_1\<,y^{}_2,y^{}_3,\dots$  on $Y$ approaching $y^{}_0$ such that for all~$i$, $X_{y_i}$ is nonempty and
 $f$ is equinormalizable at each $x\in X_{y^{}_i}$ (see \ref{fibers} or \cite[p.\,271, (I.10)]{Mn}). As shown above, this implies that 
$$
\dt(X_{y^{}_i})=n-r=\dim_{K_{y^{}_0}}\! \Big((\>\>\overline{\<\<\OX\<}/\OX)_{x_0^{}} \otimes_{\sO_{y^{}_0}} K_{y^{}_0}\Big).
$$
Since for $i\gg 0$, $\dt(X_{y^{}_i})=\dt(X_{y^{}_0})$, we
conclude that, with $\sO_{y^{}_0}\set\sO_{Y\!,\>\>y^{}_0}$ and $\sO_{x^{}_0}\set\sO_{\<\<X,\>y^{}_0}$,
$$
\dt_{\>\C}(\sO_{x_0^{}}\otimes_{\sO_{y^{}_0}}k_{y_0^{}})=
\dt(X_{y^{}_0})=\dim_{K_{y^{}_0}}\! \Big((\>\>\overline{\<\<\OX\<}/\OX)_{x^{}_0} \otimes_{\sO_{y^{}_0}} K_{y^{}_0}\Big)
=\dt_{K_{y^{}_0}}(\sO_{x_0^{}}\otimes_{\sO_{y^{}_0}}K_{y^{}_0}).
$$
In view of Proposition~\ref{cluban}, Theorem~\ref{main} gives that $f$ is equinormalizable at $x_0^{}$.
\end{proof}

\begin{ind}
Theorem~\ref{main2an} below deals with a \emph{projective} analytic map
$f\colon X\to Y$, that is, a proper~map such that there exists an
$f$-ample invertible $\OX$-module
\cite[p.\,141]{BS}.  (To avoid some trivialities we assume henceforth
that $X$ is nonempty.) The theorem concerns properties which are
 local on $Y\<$, so  no essential loss in
generality results from assuming  $Y$~connected and, with 
$\mathbb P^{\>r}$ the $r$-dimensional projective space over $\C$  for some
$r>0$, that \mbox{$f=p\smcirc i\>\>$} where 
\mbox{$i\colon X\to Y\<\times\>\> \mathbb P^{\>r}$} is a closed immersion and 
$\>p\colon Y\<\times\>\> \mathbb P^{\>r}\to Y$ is the projection (see
\cite[p.\,143]{BS}). Using Serre's ``\kern.5pt GAGA" comparison theorem (see
e.g., \cite[p.\,148, Thm.\,2.6]{BS}) one sees that any fiber 
$X_y\ (y\in Y)$ is the complex subspace associated to a closed 
subscheme~$X_y^0$
of~$\>\mathbb P^{\>r}$ (Chow's theorem) and that for any coherent
$\sO_{\<\<X_{\<\<y}}$-module\/~$\mathcal M$ the definition of the
Hilbert polynomial $\sH_y(\mathcal M)$---recalled just before
Prop.\,~\ref{hilbconst}, with $\cL\set i^*\sO_{Y\<\times\> \mathbb
P^{\>r}}(1)$---makes sense in the present context: in fact if
$\gamma\colon X_y\to X_y^0$ is the canonical ringed-space map then
$\mathcal M\cong\gamma^*\mathcal M^0$ for some coherent
$\sO_{\<\<X_y^0}$-module~$\mathcal M^0$,\vspace{-.6pt} unique up to
isomorphism, and $\sH_y(\mathcal M)$ (analytic) coincides with $\sH_y(\mathcal
M^0)$~(algebraic)---the latter\vspace{.6pt} calculated with respect to 
the algebraic $\sO(1)$ on $X_y^0\subset\mathbb P^{\>r}\<$.

In particular, when $X_y$ (and hence $X_y^0$)
is reduced and $\mathcal M^0$
is the integral closure of~$\>\sO_{\<\<X_y^0}$ then
by basic properties of the association of complex spaces\vspace{-2pt}  
to finite-type $\C$-schemes 
(clearly described in \cite[Expos\'e XII]{SGA}),
$\mathcal M\set \gamma^*\mathcal M^0$ is the integral closure of~
$\>\sO_{\<\<X_y}$. (For this fact, another argument appears in the proof 
of Lemma~\ref{H_K}.)

\end{ind}

For simplicity, we assume throughout that the complex space $Y$ is
normal and irreducible. To deduce the analytic Theorem~\ref{main2an}
from its algebraic counterpart Theorem~\ref{main2} we will need some
results on projective analytic maps analogous to ones given for
schemes in \cite{EGII} and~\cite{EGIII}, and enabled by
Grauert-Remmert's analytic version \cite[p.\,36, (5.11)]{C} of Serre's
``fundamental theorem of projective morphisms"
\cite[Thm\,(2.2.1)]{EGIII}. In our description of these results,
some details will be left implicitly to the reader.

Every  finitely-presented $\sO_Y$-algebra $\sG=\oplus_{n\ge0}\,\>\sG_n$
\cite[p.\,19-01, D\'efn.\,1]{C} determines an analytic map
$g\colon\pn\sG\to Y$ such that for each Stein open $U\subset Y$ whose
closure $\overline U$ is a Stein compactum,%
\footnote{That is, a semianalytic compact set having a neighborhood basis
of open Stein subspaces of $Y\<$, for example, the inverse image of a
small closed ball under a local embedding of $Y$ into some $\C^n$.}
there is an $s>0\>$ such that the reduced space underlying $g^{-1}U$
looks like the zeros in~$U\times\>\mathbb P^{\>s}$ of a set of homogeneous
polynomials in $\Gamma(\>\>\overline{\< U}\<,\sO_Y)[T_0,T_1,\dots,
T_s]$ generating the kernel of a surjection $(\sO_Y|_{\>\overline{\<
U\<}\>})[T_0,T_1,\dots, T_s]\twoheadrightarrow
\sG|_{\overline U}$. 
(Details are in \cite[p.\,36]{Bi}.)  For instance,
$$
\pn \sO_Y[T_0,T_1,\dots, T_s]=Y\<\times\> \mathbb P^{\>s}.
$$
If $g\colon Y_1\to Y$ is an analytic map, then there is a natural
$Y_1$-isomorphism
\begin{equation}\label{baseproj}
\pn g^*\sG\iso(\pn\sG)\times_Y Y_1.
\end{equation}
To see this, restrict to a suitable Stein open subset of $Y\<\<$, 
cover $\pn \sG$ by open subspaces of the form $\specan \sG'$ (see
\cite[p.\,19-02, D\'efn.\,2]{C}, \cite[p.\,37, top]{Bi},
\cite[(3.1.4)]{EGII}) and then using that 
$\specan g^*\sG'=\specan \sG\times_Y Y_1$
\cite[p.19-03, Prop.\,2(iv)]{C} argue by pasting, 
as in \mbox{\cite[p.\,10-05, Lemme 2.3]{C}.}

\pagebreak[3]

Set $Z\set Y\<\times\> \mathbb P^{\>r}$ and let $p\colon Z\to Y$ 
be the projection.
If $Y$ is a Stein space, then the calculations leading to
\cite[(2.1.12)]{EGA} carry over to the analytic context, giving a
natural isomorphism
$$
H^0(Y,\sO_Y)[T_0,T_1,\dots, T_r]\iso
H^0\Big(Z,\,\oplus_{n\ge0}\;\sO_Z(n)\Big),
$$
and hence, for general $Y\<$, a natural isomorphism of $\sO_Y$-algebras
$$
\sO_Y[T_0,T_1,\dots, T_r]\iso\oplus_{n\ge0}\:p_*\Big(\sO_Z(n)\Big).
$$
 
For $f\colon X\overset{i}\hookrightarrow Z\xrightarrow{p} Y$ 
as above, since $f_*\OX(n)$ is coherent, the
image~$\sG$ of the natural graded homomorphism
$\pi\colon\oplus_{n\ge0}\;p_*(\sO_Z(n))\to \oplus_{n\ge0}\;
f_*(\OX(n))$ is finitely presented, see \cite[p.\,2,
Prop.\,1.4]{TL}. For all $n\gg 0$ the graded component $\pi_n$ of
$\pi$ is \emph{surjective,} because if $I$ is the (coherent) kernel of
the natural map $\sO_Z\to i_*\OX$ then we have the natural exact
sequence
$$
0\to I(n)\to \sO_Z(n)\to i_*\OX(n)\to 0,
$$
and for $n\gg0$, $R^1\<p_*I(n)=0$ (by the above-mentioned theorem of 
Grauert-Remmert). By \cite[p.\,39, (5.13)(2)]{Bi} we have then
a $Y\<$-isomorphism $X\iso\pn\sG$. 

For any $y\in Y\<$, the fiber $X_y$ is 
$\pn \sG(y)$\vspace{.4pt} where $\sG(y)\set\sG_y\otimes_{\sO_{Y\!,\>\>y}}\C$ 
($\C$ being identified with the residue field
of~$\sO_{Y\!,\>\>y}$), see \eqref{baseproj}. 
The construction of $\pn$ in \cite[p.\,36]{Bi} shows that 
$X_y$ is the complex space associated to the projective \mbox{$\C\>$-scheme} 
$X_y^0\set\proj \sG(y)$, the closed fiber of the natural map 
$f_y\colon\proj \sG_y\to\spec {\sO_{Y\!,\>\>y}}$ ($\sG_y$ being the stalk at $y$ of $\sG$). So if $X_y$ is
reduced and $y^0$ is the
closed point of $\spec {\sO_{Y\!,\>\>y}}$ then by the above
remarks the analytic and algebraic Hilbert
polynomials of the integral closures of the structure sheaves on these
fibers coincide:
\looseness=-1
\begin{equation}\label{samehilb}
\sH_y(\>\>\overline{\<\<\sO_{\<\<X_{\<y}}\<\<}\>\>)=
\sH_{y^0}(\>\>\overline{\<\<\sO_{\<\<X^{\<0}_{\<\<y}}\<\<}\>\>).
\end{equation}

\vspace{1pt}

Assume henceforth that $f$ \emph{satisfies} $\condan$. By ~\ref{Nishimura} and \cite[p.\,184, Cor.\,(ii)]{Ma}, $X$~is \emph{reduced.}  Let $\eta=\eta_y$ be
the generic point of $\spec\sO_{Y\!,\>\>y}$, and with $K_y$ the
fraction field of\/ $\sO_{Y\!,\>\>y}$ let~$X_\eta$ be the $K_y$-scheme
$\proj(\sG_y\otimes_{\sO_{Y\!,\>\>y}}K_y)$.\vspace{-1pt}  The integral closure
$\>\overline{\<\sO_{\<\<X_\eta}\<\<}\>\>$ of $\sO_{\<\<X_\eta}$ is a
coherent $\sO_{\<\<X_\eta}$-module; and as in the paragraph preceding
Proposition~\ref{hilbconst} (with $\cL=\sO_{\<\<X_\eta}(1)$) we have the Hilbert polynomial
$\sH_\eta(\>\>\overline{\<\<\sO_{\<\<X_\eta}\<\<} \>\>)$.

\addtocounter{slem}{2}

\begin{slem}\label{Hilbeq} 
In the preceding circumstances,
$\sH_y(\>\>\overline{\<\<\OX_{\<\<y}\<\<}\>\>)=
\sH_\eta(\>\>\overline{\<\<\sO_{\<\<X_\eta}\<\<} \>\>)$\ 
iff\/ there is an open neighborhood $U$ of $y$ 
such that $f$ is equinormalizable everywhere on $f^{-1}U$.
\end{slem}

\begin{proof} We will apply Theorem~\ref{main2}. To relate the algebraic and complex-analytic setups,  note
first that the natural ringed-space map\vspace{-.4pt} 
$\varphi\colon X_y=\pn\sG(y)\to\proj\sG(y)=X_y^0$ maps the points of~$X_y$ bijectively to the closed points of $X^0_y$,\vspace{-2pt} and that for each $x\in X_y$ the induced maps are isomorphisms of completions
$\hat\varphi_x\colon\hat\sO_{\<\<X^0_y\<,\>\varphi(x)}\iso  \hat\sO_{\<\<X_y\<,\>x}$  (see \cite[Expos\'e XII, Thm.\,1.1]{SGA}). Since $\sO_{\<\<X_y\<,\>x}$ is reduced
(by $\condan$) and excellent, it follows that $X^0_y$, the closed fiber of the above map\/ $f_y\colon\proj \sG_y\to\spec {\sO_{Y\!,\>\>y}}$, is reduced.\vspace{1pt}

More generally, setting $Z_y\set \proj \sG_y$ we claim that \emph{the natural map\/ $\sO_{\<Z_y,\varphi(x)}\to \sO_{\<\<X\<,\>x}$ induces an isomorphism\vspace{-.5pt} of completions.} This granted, it follows (since $f$ satisfies $\condan$) that $\sO_{\<Z_y,\varphi(x)}$~is equidimensional\vspace{-1pt}
 and flat over~$\sO_{Y,y}$. Since the points $\varphi(x)$\vspace{.6pt} are all the closed points of~$Z_y$,\vspace{-.6pt} one concludes, using~\ref{Nishimura},
that $f_y$ \emph{satisfies} $\cond$. Then in view of~\eqref{samehilb}, Theorem~\ref{main2} gives that  $\sH_y(\>\>\overline{\<\<\OX_{\<\<y}\<\<}\>\>)=\sH_\eta(\>\>\overline{\<\<\sO_{\<\<X_\eta}\<\<} \>\>)$ \emph{iff\/ $f_y$ is equinormalizable.}  

By~\ref{simnorlocal}, if $f_y$ is equinormalizable then its completion\vspace{-1pt} 
$\hat f_{\<x}^{}\colon\spec \widehat{\sO_{\<\<X\<\<,\>x}}\to\spec \widehat{\sO_{Y\!,\>\>y}}$  is equinormalizable, for all $x\in X_y$.
Conversely, if $\hat f_{\<x}^{}$ is equinormalizable\vspace{-1.5pt} then, with $C\set\sO_{Z_y,\>\varphi(x)}$ and 
$D\set\overline C$, one has, as in the proof  of~\ref{simnorlocal}, that $\hat D$ is $\widehat{\sO_{Y\!,\>\>y}}$-flat, whence $D$ is $\sO_{Y\!,\>\>y}$-flat, and that the closed fiber of the obvious map\ $\bar f_{\<x}^{}\colon\spec D\to\spec\sO_{Y\!,\>\>y}$ is geometrically normal, whence, by~\ref{Nishimura}, all the fibers of $\bar f_{\<x}^{}$ are geometrically normal. 
So if  $\hat f_{\<x}^{}$ is equinormalizable for all $x\in X_y$ then it results from
\vspace{-.8pt}
~\ref{CorSnormal} that $f_y$ is equinormalizable. Thus \emph{$f_y$ is equinormalizable iff\/ $\hat f_{\<x}^{}$ is equinormalizable for all\/ $x\in X_y$.}\vspace{1pt}

Now \emph{$\hat f_{\<x}^{}$ is equinormalizable for all\/ $x\in X_y$ iff\/ $f$ is equinormalizable at every such\/ $x.$} Indeed, 
using~\ref{Nishimura} one finds as above that since $f$ satisfies $\condan$ therefore the induced map $\spec \sO_{\<\<X\<,\>x}\to\spec\sO_{Y\!,\>y}$ satisfies $\cond$; and as in the preceding paragraph, we see that this map is equinormalizable---i.e., $f$ is equinormalizable at~$x$---if and only if its completion at $x$, which, by the above claim, is  \emph{the same\/ $\hat f_{\<x}^{}$ as before,} is equinormalizable.\vspace{1pt} 

In conclusion, $\sH_y(\>\>\overline{\<\<\OX_{\<\<y}\<\<}\>\>)=\sH_\eta(\>\>\overline{\<\<\sO_{\<\<X_\eta}\<\<} \>\>)$ \emph{iff\/ $f$ is equinormalizable at every\/ $x\in X_y,$}
i.e., by ~\ref{equinorm local}, \emph{iff\/ $f$ is equinormalizable on some  neighborhood of\/ $f^{-1}y,$} i.e., since $f$ is proper, \emph{iff\/ $y$ has a neighborhood  on whose inverse image $f$ is equinormalizable,} as~asserted.\vspace{1pt}

As for the the claim, let $U$ be any Stein neighborhood of $y$ whose closure $\overline U$ is a Stein compactum, set $A_U\set\Gamma(U,\sO_Y)$ and  $S_U\set\Gamma(\overline U,\sG)\otimes_{\Gamma(\overline U,\sO_Y)}A_U$,
so that $f^{-1}U$ is the complex space associated to the $A_U$-scheme $X_U\set\proj S_U$, coming with a canonical ringed-space map $\varphi^{}_U\colon f^{-1}U\to X_U$
\cite[p.\,36]{Bi}.
As $\sG_y=\dirlm{U} S_U$,\vspace{.6pt} it follows for any~$x\in X_y$ that \mbox{$\sO_{\<Z_y,\varphi(x)}=\dirlm{U}\sO_{\<\<X_U\<,\>\varphi^{}_U(x)}$.}\vspace{2pt} With $\fm_x$ (resp.~$\fm_U$) the maximal ideal of $\sO_{\<\<X,\>x}$ (resp~$\sO_{\<Z_U\<,\>\varphi^{}_U(x)})$,  $\varphi^{}_U$ induces for any $n>0$\vspace{1pt}
 an isomorphism  
$\sO_{\<\<X_U\<,\>\varphi^{}_U(x)}/\fm_U^n\iso \sO_{\<\<X,\>x}/\fm_x^n$\vspace{1pt}
 \cite[p.\,1, (1.1)]{Bi}; taking the direct limit over $U$ gives, with $\fm$ the maximal ideal of $\sO_{\<Z_y,\varphi(x)}$,  an isomorphism
$\sO_{\<Z_y,\varphi(x)}/\fm^n\iso \sO_{\<\<X,x}/\fm_x^n$, whence
the claim.
\end{proof}

\begin{slem}\label{H_K}
 The Hilbert polynomial 
$\sH_\eta(\>\>\overline{\<\<\sO_{\<\<X_\eta}\<\<} \>\>)$ does not depend on\/~$y$. \end{slem}

\begin{proof}

Let $\>\overline{\<\OX\<\<}\>\>$ be the integral closure of $\OX$, a coherent $\OX$-module \cite[p.\,21-09, Cor.\,2]{C}. For each $n>0$ the $\sO_Y$-module 
$\overline\sG_n\set f_*(\>\>\overline{\<\<\OX\<\<}\>\>(n))$
is coherent, hence---as in the proof of Theorem~\ref{mainan}---locally free on a dense open subset of~$Y\<$, of rank, say $r_n\>$;
and \emph{for each\/ $y\in Y,$ and\/ $K_y$ the fraction field of\/ $\sO_{Y\!,\>\>y},$ the dimension of the\/ $K_y$-vector space\/ 
$(\overline\sG_n)_y\otimes_{\sO_y} K_y$ is\/ $r_n$}.

Let us describe the (algebraically defined) polynomial
$\sH_\eta(\>\>\overline{\<\<\sO_{\<\<X_\eta}\<\<} \>\>)$ in terms of
the (analytically defined) integers $r_n$---which do not depend on~$y$. 
Let  $U\subset Y$ be any Stein open set whose closure $\overline U$ is a Stein compactum, and let $A_U$, and $\varphi^{}_U\colon f^{-1}U\to X_U$ be as in the preceding proof.
The canonical very ample invertible sheaves on $X_U$ and~$f^{-1}U$ are related in the obvious way: $\varphi_U^*\sO_{\<\<X_U}(1)=\sO_{\<\<X_U}(1)$. It follows trivially from the definition \cite[p.\,1, Satz (1.1)]{Bi} that the $U$-analytic space
associated to $\spec A_U$ is $U$ itself, together with the canonical
ringed-space map $i^{}_U\colon U\to\spec A_U$.\vspace{1pt}

\enlargethispage*{5pt}
Let $\>\overline{\<\<\<\sO_{\<\<X_U}}\>\>$ be the integral 
closure of $\sO_{\<\<X_U}$.
Then $\varphi_U^* \,\overline{\<\sO_{\<\<X_U}\<\<}\>\>=\overline{\sO_{f^{-1}U}}$: 
indeed, for $x\in f^{-1}U$ the local rings $\sO_{\<\<X_U\<,\>\>\varphi^{}_U(x)}$ and $\sO_{\<\<X\<,\>x}$ are excellent \cite[p.\,153]{Bi2}, of equicharacteristic zero, and have the same completion \cite[p.\,1, Satz (1.1)]{Bi}, so the canonical
map from the first  to the second
is \emph{regular} (flat, with geometrically regular fibers)
\cite[(6.6.1)]{EGA}; and \cite[(6.14.1)]{EGA} implies the assertion.
(One could also argue on the basis of various ``comparison'' results
in \cite[\S\S2--3]{Bi}.)  Hence by \cite[p.\,16, (4.2)]{Bi} \vspace{.6pt}
(``relative GAGA''),
if $f_U\colon X_U\to\spec A_U$ is the projection then for
all $n\ge 0$,
$
i_U^*f_{U*}\Big(\>\overline{\<\<\<\sO_{\<\<X_U}}\>(n)\Big)=
f_*\Big(\>\>\overline{\<\<\sO_{f^{-1}U}\<\<}\>\>(n)\Big).
$
\pagebreak[3]

Now $f_{U*}\Big(\>\overline{\<\<\<\sO_{\<\<X_U}}\>(n)\Big)$ is a coherent sheaf
on $\spec A_U$, associated to the finitely-generated $A_U$-module $\Gamma\Big(X_U, \,\overline{\<\<\sO_{\<\<X_U}\<\<}\>(n)\Big)$,
so its restriction to some nonempty open subscheme is locally free, necessarily of constant rank $r_n$, since that is the rank of the restriction of $f_*\Big(\>\>\overline{\<\<\sO_{f^{-1}U}\<\<}\>\>(n)\Big)$ to some dense open subset of $U$. Thus  \emph{$r_n$~is the rank of the finitely generated
$A_U$-module} $\Gamma\Big(X_U, \,\overline{\<\<\sO_{\<\<X_U}\<\<}\>(n)\Big)$.

For any $y\in U$ this last statement remains true under base change
from $A_U$ to $A_y$, where $A_y$ is the localization of $A_U$ at the
maximal ideal~$i_U(y)$ of~$A_U$ (because $\Gamma$ ``commutes'' with flat base change \cite[p.\,366, (9.3.3)]{EGA1}, and integral closure is compatible with localization, or more generally, with flat base change having  geometrically normal fibers \cite[(6.14.4)]{EGA}). Furthermore, as before, the natural local homomorphism 
$A_y\to\sO_{Y\!,\>\>y}$ is regular; and therefore the statement remains
true under the base change from $A_y$ to $\sO_{Y\!,\>\>y}$. 
Now, fixing $y\in Y$, let $V$ be a neighborhood such that the finitely presented 
$\sO_{Y\!,\>\>y}$-algebra $\sG_y$ is specified by generators and relations
which are stalks of generators and relations of the $A_V$-algebra $\Gamma(V, \sG)$,
and choose a neighborhood~$U$ as above such that $\overline U\subset V$.
Then  $X_U\otimes_{A_U}\sO_{Y\!,\>\>y}$ is the above scheme $X_y$, and  
$r_n$ \emph{is the rank of the\/  $\sO_{Y\!,\>\>y}$-module} 
$\Gamma\Big(X_y, \,\overline{\<\<\sO_{\<\<X_y}\<\<}\>(n)\Big)$. Hence
for  all  $n\gg0$,
$$
\sH_\eta(\>\>\overline{\<\<\sO_{\<\<X_\eta}\<\<} \>\>)(n)=\dim_{K_y}\Gamma\Big(X_\eta, \,
 \overline{\<\<\sO_{\<\<X_\eta}\<\<}\>(n)\Big)=r_n,
$$
which, again, does not depend on $y$.
\end{proof}

\begin{thm}\label{main2an}
Let\/ $f\colon X\hookrightarrow Y\<\times \mathbb P^{\>r}\twoheadrightarrow Y$ be a projective analytic map 
satisfying\/ $\condan\>,$ with\/ $Y$ irreducible. Then\/ $f$ is equinormalizable if and only if\/ $\sH_y(\>\>\overline{\<\<\sO_{\<\<X_y}\<\<}\>\>)$ is the same for all\/ $y\in Y\<$.
\end{thm}

\begin{proof} If $f$ is equinormalizable and $\nu\colon\BX\to X$ is a normalization then $f\nu$ is flat over~$Y\<$, and hence $\overline{\OX\<\<}\>\>=\nu_*\sO_{\<\BX}$ is $Y$-flat (its stalk at $x\in X$ being the product of the stalks of $\sO_{\<\BX}$ at the finitely many points in $\nu^{-1}x$). So by a theorem of Grauert
\cite[p.\,134, Thm.\,4.12(iii)]{BS}, $\sH_y\Big((\>\>\overline{\<\<\OX\<\<}\>\>)_y\Big)$ is independent of $y$. Furthermore, the fibers of $f\nu$ normalize those of $f$, and hence one gets a natural isomorphism $(\>\>\overline{\<\<\OX\<\<}\>\>)_y\iso\>\>\overline{\<\<\sO_{\<\<X_y}\<\<}\>\>$. Thus 
$\sH_y(\>\>\overline{\<\<\sO_{\<\<X_y}\<\<}\>\>)$
 is independent of~$y$.
 
 Suppose conversely that $\sH_y(\>\>\overline{\<\<\sO_{\<\<X_y}\<\<}\>\>)$ does not depend on $y$. By  Prop.\,\ref{fibers} or \cite[p.\,271, (I.10)]{Mn}, there is a 
$y^{}_1\in Y$ such that 
 $f$ is equinormalizable at each $x\in X_{y^{}_1\<}$. For any $y^{}_0\in Y$
 Lemmas~\ref{Hilbeq} and ~\ref{H_K} give
 $$
 \sH_{y^{}_0}(\>\>\overline{\<\<\sO_{\<\<X_{y^{}_0}}\<\<}\>\>)=
 \sH_{y^{}_1}(\>\>\overline{\<\<\sO_{\<\<X_{y^{}_1}}\<\<}\>\>)=
 \sH_{\eta^{}_1}(\>\>\overline{\<\<\sO_{\<\<X_{\eta^{}_1}}\<\<}\>\>)=
 \sH_{\eta^{}_0}(\>\>\overline{\<\<\sO_{\<\<X_{\eta^{}_0}}\<\<}\>\>),
$$ 
and then Lemma~\ref{Hilbeq} shows that $f$ is equinormalizable over 
a neighborhood of~${y^{}_0}$.
\end{proof} 
\newpage

\setcounter{footnote}{0}
\renewcommand{\thefootnote}{\fnsymbol{footnote}}


\begin{thebibliography}{EGA\,I}


\bibitem[BG]{BG} C.\,Br\"ucker, G.-M.\,Greuel, 
\emph{Deformationen isolierter Kurvensingularitäten mit eingebetteten Komponenten.} Manuscripta Math. {\bf70} (1990),  93--114. 

\bibitem[BGG]{BGG} J.\,Brian\c con, A.\,Galligo, M.\,Granger, 
\emph{D\'eformations \'equisinguli\`eres des germes de courbes gauches 
r\'eduites.} 
M\'em.~Soc.~Math.~France (N.S.) 1980/81, no. 1.

\bibitem[Bi]{Bi} J.\, Bingener,
\emph{Schemata \"uber Steinschen Algebren.}  Schr.\ Math.\ Inst.\
Univ.\ M\"unster (2) Heft 10, 1976.



\bibitem[Bi2]{Bi2} \bysame,
\emph{Holomorph-pr\" avollst\" andige Restr\" aume zu
analytischen Mengen in Steinschen R\" aumen.}  J.~Reine
Angew.~Math. {\bf285} (1976), 149--171. 



\bibitem[BKKN]{BKKN} R.\, Berger, R.\,Kiehl, E.\,Kunz, H.-J.\,Nastold,
 \emph{Differentialrechnung in der  analytischen Geometrie.}
 Lecture Notes in Math., no.~38,
 Springer-Verlag, Berlin Heidelberg New York, 1967.
 


\bibitem[Bo]{Bo} N.\, Bourbaki,
 \emph{Alg\`ebre Commutative, Chapitre 3.}
 Actuali\'es Sci.~et~Ind., no.~1293,
 Herrmann, Paris, 1961.
 
 \bibitem[BR]{BR} A.\,Brezuleanu, C.\, Rotthaus,
 \emph{Eine Bemerkung \"uber Ringe mit geometrisch normalen formalen Fasern.}
 Arch. Math. {\bf 39} (1982), 19-27.

\bibitem[BS]{BS} C.\,B\u anic\u a, O.\,St\u an\u a\c sil\u a, 
\emph{Algebraic methods in the global theory of complex spaces.} Editura Academiei, Bucharest; John Wiley \& Sons, London-New York-Sydney, 1976.

\bibitem[C]{C} S\'eminaire Henri Cartan, 13i\`eme ann\'ee: 1960/61. 
\emph{Familles d'espaces complexes et fondements de la g\'eom\'etrie 
analytique.} Fasc. 1 et 2: Exp. 1--21.  
2i\`eme \'edition, corrig\'ee. 
\'Ecole Normale Sup\'erieure Secr\'etariat math\'ematique, Paris, 1962 
 
\bibitem[Fi]{Fi} G.\,Fischer,
 \emph{Complex Analytic Geometry.}
 Lecture Notes in Math., no.~538,
 Springer-Verlag, Berlin Heidelberg New York,
 1976.

\bibitem[Fr]{Fr} J.\,Frisch,
 \emph{Points de platitude d'un morphisme
 d'espaces analytiques complexes.}
 Invent. Math. {\bf 4} (1967), 118--138.

\bibitem[G]{G} A.\,Grothendieck,
\emph{Techniques de construction et th\'eor\`emes d'existence 
en g\'eom\'etrie alg\'ebrique. IV. Les sch\'emas de Hilbert}.
 S\'eminaire Bourbaki. Vol.~6.  
Ann\'ee 1960/61. Expos\'es 205--222. 
Soci\'et\'e Math\'ematique de France, Paris, 1995. iv+286 pp.


\bibitem[EGI]{EGA1} A.\,Grothendieck, J.\,Dieudonn\'e,
 \emph{\'El\'ements de  G\'eom\'etrie Alg\'ebrique I.
 }
 Springer-Verlag,  New York,  1971.

\bibitem[EGII]{EGII} \bysame,
 \emph{\'El\'ements de  G\'eom\'etrie Alg\'ebrique III,
 \'Etude globale  \'el\'ementaire de quelques classes de morphismes.}
 Publications Math. I.H.E.S. {\bf 8} (1961).\footnote[2]{Downloadable from 
{\tt <http://www.numdam.org/numdam-bin/feuilleter?j=PMIHES\&sl=0>}.}



\bibitem[EGIII]{EGIII} \bysame,
 \emph{\'El\'ements de  G\'eom\'etrie Alg\'ebrique III,
 \'Etude cohomologique des faisceaux coh\'erents.}
 Publications Math. I.H.E.S. {\bf 11} (1961), {\bf17} (1963).\footnotemark[2]

\bibitem[EGA]{EGA} \bysame,
 \emph{\'El\'ements de  G\'eom\'etrie Alg\'ebrique IV,
 \'Etude locale des sch\'emas et des morphismes de sch\'emas.}
 Publications Math. I.H.E.S. {\bf 24} (1965), {\bf 28} (1966),
 {\bf 32} (1967).\footnotemark[2]

\bibitem[SGA]{SGA} \bysame, Mich\`ele Raynaud,
\emph{Rev\^etements \'etales et groupe fondamental} (SGA 1). S\'eminaire de g\'eom\'etrie alg\'ebrique du Bois Marie 1960--61. Updated reprint of Lecture Notes in Math., 224, Springer, Berlin. Documents Math\'ematiques (Paris), 3. Soc.\,Math\, de France, Paris, 2003.

\bibitem[GrS]{GS}\bysame, H.\,Seydi,
\emph{Platitude d'une adh\'erence sch\'ematique et lemme de Hironaka g\'en\'eralis\'e.}  
Manuscripta Math. {\bf 5} (1971), 323--339.



\bibitem[GR]{Gr} H.\,Grauert and R.\,Remmert,
 \emph{Coherent Analytic Sheaves.}
 Springer-Verlag,  New York,  1984.


\bibitem[Ki]{Kl} R.\,Kiehl,
 \emph{Ausgezeichnete Ringe in der nichtarchimedischen analytischen Geometrie.}
  J. Reine Angew. Math. {\bf  234}  (1969), 89--98.

\bibitem[Ku]{Kz} E.\,Kunz,
 \emph{On noetherian rings of characteristic $p$.}
  Amer. J.  Math. {\bf  98}  (1976), 999--1013.

\bibitem[Lp]{Lp} J.\,Lipman,
\emph{Free derivation modules on algebraic varieties.}
 Amer.~ J.~ Math. {\bf 87} (1965), 874--898.

\bibitem[Ma1]{M1} H.\,Matsumura,
 \emph{Commutative Algebra, second edition.}
 W.A. Benjamin, New York, 1980.

\bibitem[Ma2]{Ma} \bysame,
 \emph{Commutative Ring Theory.}
 Cambridge Studies in Adv. Math. {\bf 8}, Cambridge University Press, 1989.
 
\bibitem[Mn]{Mn} M.\,Manaresi,
\emph{Sard and Bertini type theorems for complex spaces.}
  Ann. Mat. Pura Appl. (4) {\bf 131} (1982), 265--279. 

 \bibitem[Mr]{Mr} J.\,Marot,
 \emph{Sur les anneaux universellement japonais.}
 Bull.\,Soc.\,Math.\,France {\bf 103} (1975), 103--111.

\bibitem[Ng]{Nag} M.\,Nagata,
 \emph{Local Rings.}
 Wiley, New York, 1962.

\bibitem[Ni]{Ni} J.\,Nishimura,
 \emph{On ideal-adic completion of noetherian rings.}~
 J.~Math.~Kyoto~Univ.~{\bf 21}~(1981), \mbox{153--169.}

\bibitem[R1]{R1} L.\,J.\,Ratliff,
 \emph{On quasi-unmixed local domains, the altitude formula, and the 
chain condition for prime ideals \textup{(I)}.}    
 Amer.~J.~Math. {\bf 91} (1969), 508--528. 

\bibitem[R2]{R2} \bysame,
 \emph{On quasi-unmixed local domains, the altitude formula, and the 
chain condition for prime ideals \textup{(II)}.}    
 Amer.~J.~Math. {\bf 92} (1970), 99--144.


\bibitem[R3]{R3} \bysame,
 \emph{Locally quasi-unmixed noetherian rings and ideals of the principal class.}
 Pacific J.~Math. {\bf 52} (1974), 185--205.


\bibitem[T1]{T1} B.\,Teissier,
 \emph{The hunting of invariants in the geometry of discriminants.}
 Real and Complex Singularities (Proc. Ninth Nordic Summer
 School/NAVF Sympos. Math., Oslo, 1976), pp.\:565--678.

\bibitem[T2]{T2} \bysame,
 \emph{Resolution simultan\'ee  I, II.}
 S\'eminaire sur les Singularit\'es des Surfaces.
 Lecture Notes in Math., no.~777,
 Springer, Berlin, 1980. 

\bibitem[TL]{TL} \bysame, M.\, Lejeune-Jalabert,
\emph{Contribution \`a l'\'etude des singularit\'es du point de vue du polygone de Newton.} Th\`ese d'\'Etat, Universit\'e de Paris VII, Paris, 1973.

\bibitem[Wa]{Wa} J.\,Wahl,
\emph{Equisingular deformations of normal surface singularities.}   Ann. of Math. (2) {\bf 104} (1976),  325--356.

\bibitem[Z1]{Z1} O.\,Zariski,
\emph{Studies in equisingularity. I. 
Equivalent singularities of plane algebroid curves.}  
Amer.~J.~ Math. {\bf  87}  (1965), 507--536.


\bibitem[Z2]{Z2} \bysame,
\emph{Studies in equisingularity. III. 
Saturation of local rings and equisingularity.}  
Amer.\ J.\ Math. {\bf  90}  (1968), 961--1023.


\end{thebibliography}
\end{document}